\documentclass{amsart}
\usepackage{amssymb}
\usepackage{graphicx}
\usepackage{amscd}
\usepackage{amsmath}

\newtheorem{theorem}{Theorem}
\theoremstyle{plain}

\newtheorem{corollary}{Corollary}

\newtheorem{definition}{Definition}

\newtheorem{lemma}{Lemma}

\newtheorem{proposition}{Proposition}
\newtheorem{remark}{Remark}

\numberwithin{equation}{section}

\begin{document}
\title[Koszul algebras and sheaves]{Koszul algebras and Sheaves over \\
Projective Space}
\author{Roberto Mart\'{\i}nez-Villa}
\address{Instituto de Matem\'{a}ticas de la UNAM, Unidad Morelia, Apartado
Postal 61-3, 58089 Morelia, Mich., M\'{e}xico}
\email{mvilla@matmor.unam.mx}
\thanks{Part of this paper was written during my visit to Syracuse
University on August, 2002, I thank Dan Zacharia for his kind hospitality,
for many useful discussions and his suggestions to improve this paper. I
thank CONACYT for funding the research project.}
\date{August 21, 2003}
\subjclass{14, 16, 18}

\begin{abstract}
We are going to show that the sheafication of graded Koszul modules  $%
K_{\Gamma }$ over $\Gamma _{n}=K\left[ x_{0},x_{1}...x_{n}\right]
$ form an important subcategory $\overset{\wedge }{K}_{\Gamma }$
of the coherents sheaves on projective space, $Coh(P^{n}).$ One
reason is that any coherent sheave over $P^{n}$ belongs to
$\overset{\wedge }{K}_{\Gamma }$up to shift.

More importantly, the category $K_{\Gamma }$ allows a concept of almost
split sequence obtained by exploiting Koszul duality between graded Koszul
modules over $\Gamma $ and over the exterior algebra $\Lambda .$ This is
then used to develop a kind of relative Auslander-Reiten theory for the
category $\mathit{Coh(P}^{n})$, with respect to this theory, all but
finitely many Auslander-Reiten components for $\mathit{Coh(P}^{n})$ have the
shape \textit{ZA}$_{\infty }.$ We also describe the remaining ones.
\end{abstract}

\maketitle

\section{\protect\bigskip Introduction}

The aim of the paper is to apply the methods of non commutative algebra, and
very particularly, of finite dimensional algebras to the study of sheaves on
projective space as well as to some non commutative generalizations. The
starting point for this approach was the paper by Bernsntein Gelfand Gelfand
in which they proved that there exists a derived equivalence of categories: $%
\underline{\mathit{mod}}
$ $_{\Lambda }\cong \mathit{D}^{b}(Coh(P^{n}),\left[
BGG\right] $ between the stable category $\underline{\mathit{mod}}$ $_{\Lambda
}$ of the exterior algebra $\Lambda =K<x_{0},x_{1},...x_{n}>/<x_{i}^{2},$ $%
x_{i}x_{j}+x_{j}x_{i}>$ and the derived category of bounded complexes of
coherent sheaves over projective space. It is of interest to know which
modules over the exterior algebra correspond to particular sheaves over
projective space, for example: which modules correspond to locally free
sheaves? which to torsion free sheaves?

The answer to the first question was given by J. Bernstein and S. Gelfand $%
\left[ BG\right] $ who proved that vector bundles over $P^{n}$ correspond to
"nice" modules over the exterior algebra. We will give here a different
characterization of the modules over the exterior algebra corresponding to
the locally free sheaves on projective space.

\section{Some basic facts on vector bundles over projective space}

For the convenience of the reader, and to fix notation, we recall in this
section standard theorems of algebraic geometry on sheaves and vector
bundles over projective space.

We will consider graded quiver algebras $\Gamma =KQ/I,$ where $K$ is a
field, $Q$ a finite connected quiver and $I$ an homogeneous ideal in the
grading given by path length, $I\subset J^{2},$ where $J$ is the ideal
generated by the arrows. Let $Gr\Gamma $ be the category of graded modules
and degree zero maps. In order to prove the main results it seems to be
essential to assume $\Gamma $ is noetherian hence; we will assume it. We
will denote by $gr\Gamma $ the full subcategory of all finitely generated
graded $\Gamma -$modules. We will say that a graded modules $M$ is torsion
if it is a sum of submodules of finite length, in particular, if $M$ is
finitely generated it is torsion when it is of finite length. For a given
module $M$ we will denote by $t(M)$ the sum of all submodules of finite
length.

Given two graded $\Gamma-$ modules $X$ and $Y$ we denote by $Hom_{\Gamma
}(X,Y)_{n}$ the set of all maps $f:X\rightarrow Y$ such that $%
f(X_{i})\subset Y_{i+n}$ by $Hom_{Gr\Gamma}(X,Y)$ the maps in degree zero
and $Hom_{\Gamma }(X,Y)=\underset{i\in Z}{\oplus}Hom_{\Gamma}(X,Y)_{i}$

Let $QGr\Gamma$ be the quotient category with the same objects as $Gr\Gamma,$
denote by $\pi:$ $Gr\Gamma\rightarrow QGr\Gamma$ the quotient functor, the
maps in $QGr\Gamma$ are defined as $Hom_{QGr\Gamma}(\pi(X),\pi(Y))_{0}=%
\underrightarrow{\lim}Hom_{\Gamma}(X^{\prime},Y/Y^{\prime})_{0}$ where the
limit is taken over all pairs $(X^{\prime},Y^{\prime})$ of submodules of $X$
and $Y,$ respectively, such that $X/X^{\prime}$ and $Y^{\prime}$ are
torsion. We define $Hom_{QGr\Gamma}(\pi(X),\pi(Y))=\underset{n\in Z}{\oplus }%
\underrightarrow{\lim}Hom_{\Gamma}(X^{\prime},Y/Y^{\prime})_{n}.$

Given a graded module $X$ we define the truncated module $X_{\geq n}$ as
follows: $(X_{\geq n})_{j}=0$ if $j<n$ and $(X_{\geq n})_{j}=X_{j}$ if $%
j\geq n.$ The following result is well known $\left[ S2\right] ,\left[ M2%
\right] :$

\begin{proposition}
Let $\Gamma =KQ/I$ be a noetherian graded quiver algebra. Then if $X$ is a
graded finitely generated module and $Y$ is an arbitrary module, then we
have isomorphisms: $Hom_{QGr\Gamma }(\pi (X),\pi (Y))_{m}=\underrightarrow{%
\lim }Hom_{\Gamma }(X_{\geq n},Y)_{m}$ for all $m$ . $\blacksquare $
\end{proposition}

It is well known $\left[ P\right] ,\left[ S2\right] ,$ $QGr\Gamma$ has
injective envelopes and the extension groups $Ext_{QGr\Gamma}^{k}(\pi
(X),\pi(Y))_{m}$ are obtained as derived functors of $Hom_{QGr\Gamma}(%
\pi(X).\pi(Y))_{m},$ it follows: $Ext_{QGr\Gamma}^{k}(\pi(X),\pi(Y))_{m}\cong%
\underrightarrow{\lim}Ext_{\Gamma}^{k}(X_{\geq n},Y)_{m}.$

\begin{theorem}
(Serre $\left[ Se\right] $) Let $k$ be a field, $\Gamma =K\left[
x_{0},x_{1},...x_{r}\right] /I$ be the quotient of the polynomial algebra
module an homogeneous ideal $I.$ Let $X$ be the subscheme of \textsl{P}$_{r}$
defined by $I$ and \textit{O}$_{X}$ the sheaf of regular functions. Let $%
\mathit{Coh(X)}$\textit{\ be the category of coherent sheaves, O}$_{X}(n)$
the $nth-$ power of the twisting sheaf. Define a functor. $\Gamma _{\ast }:$
\textit{Coh(X)}$\rightarrow Qgr\Gamma $ by $\Gamma _{\ast }(\mathit{F)=}%
\underset{d}{\mathit{\oplus }}H^{0}(X,\mathit{F\otimes O}_{X}(d)),$ where $%
H^{0}(X,\mathit{F)}$ denotes the global sections of the sheaf $\mathit{F.}$
The functor is an equivalence of categories with inverse: $\thicksim :$ $%
Qgr\Gamma \rightarrow \mathit{Coh(X)}$ defined as follows:

The set $\left\{ X_{f}\right\} _{f\in \Gamma }$ with $X_{f}=\left\{ x\in
X\mid f(x)\neq 0\right\} $ forms a basis of open sets for the topology of $%
X. $ The structure sheaf \textit{O is defined by O}$(X_{f})=\Gamma \left[
f^{-1}\right] _{0},$ the degree zero part of $\Gamma \left[ f^{-1}\right] ,$
which consists of all rational functions $g/h$ in $K(X)$ having no pole on $%
X_{f}$, the polynomials $g,h$ homogeneous with $deg(g)=deg(h).$ If $M$ is a
graded module define the sheaf $\overset{\thicksim }{M}$ by $\overset{%
\thicksim }{M}(X_{f})=(\Gamma \left[ f^{-1}\right] \otimes M)_{0}=M\left[
f^{-1}\right] _{0}.\blacksquare $
\end{theorem}

For a coherent sheaf $\mathit{F}$ cohomology groups $H^{q}(\mathit{F)}$ were
defined by $\overset{\vee}{C}$ech, we refer to the reader to Serre%
\'{}%
s paper for definitions and for the proof of the isomorphisms: $H^{q}(%
\mathit{F)\cong} H^{q}(X,\mathit{F)}$, where $H^{q}(X,\mathit{-)}$ denotes
the derived functor of the global sections functor: $H^{0}(X,\mathit{-)}$ .

The groups $H^{q}(X,\mathit{F)}$ can be determined in terms of the extension
groups as follows:

Let $M$ be a graded $\Gamma-$ module such that $\overset{\thicksim}{M}=%
\mathit{F.}$ It was proved in $\left[ Se\right] $ , $H^{q}(X,\mathit{F)=}$ $%
Ext_{QGr\Gamma}^{q}(\mathit{O}_{X},\mathit{F})_{0},$ the degree zero part of
the extension group.

Let $J=$ $(x_{0},x_{1},...x_{r})/I.$ Then $Ext_{QGr\Gamma}^{q}(\mathit{O}%
_{X},\mathit{F})_{0}\cong\underrightarrow{\lim}Ext_{%
\Gamma}^{q}(J^{k},M)_{0}. $ Define \underline{$H$}$^{q}(X,\mathit{F)}$ $=%
\underset{n\in Z}{\oplus}H^{q}(X,\mathit{F(n)),}$ where $\mathit{%
F(n)=F\otimes O}_{X}(n)$ and $\overset{\thicksim}{\mathit{F(n)}}\mathit{=}%
\overset{\thicksim}{\mathit{M}\left[ n\right] }$ the $nth-$ shift of $M.$

With these definitions \underline{$H$}$^{q}(X,\mathit{F)}$ $=\underset{n\in Z%
}{\oplus}\underrightarrow{\lim}Ext_{\Gamma}^{q}(J^{k},M\left[ n\right] )_{0}=%
\underrightarrow{\lim}Ext_{\Gamma}^{q}(J^{k},M).$ In particular, $%
\Gamma_{\ast}(\mathit{F)=}$\underline{$H$}$^{0}(X,\mathit{F)}$ $%
=Hom_{QGr\Gamma}(\pi(\Gamma),\pi(M)).$

Using these isomorphisms Serre duality $\left[ H2\right] $ has the following
form:

\begin{theorem}
$.$ Let $\Gamma =K\left[ x_{0},x_{1},...x_{r}\right] /I$ be the quotient of
the polynomial algebra by an homogeneous ideal and $M$ a finitely generated
graded module. Then there exists a natural isomorphism:

$Hom_{K}(\underrightarrow{\lim }Ext_{\Gamma }^{q}(J^{k},M\left[ -n\right]
)_{0},K)\cong \underrightarrow{\lim }Ext_{\Gamma }^{r-q}(M\left[ -n\right]
_{\geq k},\Gamma \left[ -(r+1)\right] )_{0}.\blacksquare $
\end{theorem}

\begin{definition}
A sheaf $\mathit{F}$ is locally free if and only if $X$ can be covered by
open sets $\mathit{U}$ for which $\mathit{F\mid }_{U}$ is a free $\mathit{O}%
_{X}\mid _{U}-$ module.
\end{definition}

We will need the following well known characterization of algebraic vector
bundles:

\begin{proposition}
$\left[ L\right] $ The functor which associates the sheaf of modules of
regular sections to a vector bundle $E$ is an equivalence of categories
between the category of algebraic vector bundles over $X$ and the category
of locally free sheaves of finite rank on $X.\blacksquare $
\end{proposition}

The following result will be also needed:

\begin{theorem}
(Serre $\left[ Se\right] )$ Let $X={P}_{r}$ be projective
space.
Then a coherent sheaf $\mathit{F}$ is locally free if and only if $H^{q}(X,$%
\textit{\ }$\mathit{F(-n))}$ $=0$ for $n>>0$ and $0\leq q<r.\blacksquare $
\end{theorem}

We can prove now the following:

\begin{proposition}
Let $\Gamma =K\left[ x_{0},x_{1},...x_{r}\right] $ be the polynomial algebra
and $\mathit{F}=\pi (M)=\overset{\thicksim }{M}$ the sheaf associated to a
finitely generated module $M$. Then $\mathit{F}$ is locally free if and only
if $\underrightarrow{\lim }Ext_{\Gamma }^{j}(M_{\geq t},\Gamma \left[ s%
\right] )_{0}=0$ for $s>>0$ and $0<j\leq r.$
\end{proposition}

\textit{Proof. }It is easy to see $M\left[ -n\right] _{\geq k}=M_{\geq k-n}%
\left[ -n\right] .$

By Serre duality:

\begin{center}
$Hom_{K}(\underrightarrow{\lim }Ext_{\Gamma }^{q}(J^{k},M\left[ -n\right]
)_{0},K)\cong \underrightarrow{\lim }Ext_{\Gamma }^{r-q}(M\left[ -n\right]
_{\geq k},\Gamma \left[ -(r+1)\right] )_{0}$

$\cong \underrightarrow{\lim }Ext_{\Gamma }^{j}(M_{\geq k-n}\left[ -n\right]
,\Gamma \left[ -(r+1)\right] )_{0}\cong \underrightarrow{\lim }Ext_{\Gamma
}^{j}(M_{\geq k-n},\Gamma \left[ n-(r+1)\right] )_{0}.$
\end{center}

Now it follows $H^{q}(X,$\textit{\ }$\mathit{F(-n))}$ $=0$ for $n>>0$ and $%
0\leq q<r$, if and only if $\underrightarrow{\lim }Ext_{\Gamma }^{j}(M_{\geq
t},\Gamma \left[ s\right] )_{0}=0$ for $0<j\leq r$ and $s>>0.$ $\blacksquare
$

\section{Koszul Algebras.}

In this section we recall some basic facts on Koszul algebras, the results
given here were proved in $\left[ BGS\right] $, $\left[ ADL\right] $, $\left[
GM1\right] $, $\left[ GM2\right] $. See also the references given there.

\begin{definition}
Let $\Lambda =KQ/I$ be a graded quiver algebra. A finitely generated module $%
M$ with minimal graded projective resolution: $...\rightarrow
P_{k}\rightarrow P_{k-1}\rightarrow ...\rightarrow P_{1}\rightarrow
P_{0}\rightarrow M\rightarrow 0$ with each $P_{j}$ finitely generated is
called quasi Koszul if and only if $J\Omega ^{k}(M)=\Omega ^{k}M\cap
J^{2}P_{k-1}$ for any $k\geq 1.$

The module $M$ is called weakly Koszul if and only if $J^{j}\Omega
^{k}(M)=\Omega ^{k}M\cap J^{j+1}P_{k-1}$ for any $j\geq 0$ and any $k\geq 1$.

The module $M$ is Koszul if and only if each $P_{j}$ is finitely generated
in degree $j.$

A graded algebra is called Koszul if and only if all graded simple
(generated in degree zero) are Koszul.
\end{definition}

\begin{definition}
Let $\Lambda $ be a graded quiver algebra and $M$ a module with minimal
graded injective co-resolution: $0\rightarrow M\rightarrow I_{0}\rightarrow
I_{1}\rightarrow ...I_{k}\rightarrow ...$, with each $I_{k}$ finitely
cogenerated.

1) We say that $M$ is quasi co-Koszul if for any $j>0$ there exist an
epimorphism:

$soc^{2}I_{j-1}/socI_{j-2}\twoheadrightarrow soc\Omega ^{-j}(M)\rightarrow
0. $

\textit{2) If for any pair of positive integers }$j$\textit{\ and }$k$%
\textit{\ there are epimorphisms}:

$soc^{k+2}I_{j-1}/soc^{k+1}I_{j-2}\twoheadrightarrow soc^{k+1}\Omega
^{-j}(M)/soc^{k}\Omega ^{-j}(M)\rightarrow 0.$

\textit{Then we say} $M$\textit{\ is weakly co-Koszul.}

3)\textit{\ If for all non negative integers }$j$\textit{\ the injective
module }$I_{j}$\textit{\ is cogenerated in degree }$j,$\textit{\ then we say
}$M$\textit{\ is co-Koszul.}
\end{definition}

\begin{remark}
Weakly Koszul (co-Koszul) modules were called strongly quasi Koszul (quasi
co-Koszul) in $\left[ GM1\right] $, $\left[ GM2\right] .$
\end{remark}

We will denote by $K_{\Lambda}$ the category of Koszul modules and degree
zero maps, by $\Gamma$ the Yoneda algebra, $\Gamma=\underset{k\geq0}{\oplus }%
Ext_{\Lambda}^{k}(\Lambda_{0},\Lambda_{0})$ and by $F:Gr_{\Lambda}%
\rightarrow Gr_{\Gamma^{op}}$ the ext functor: $F(M)=\underset{k\geq0}{\oplus%
}Ext_{\Lambda}^{k}(M,\Lambda_{0}).$

\begin{theorem}
Let $\Lambda =KQ/I$ be a Koszul algebra with Yoneda algebra $\Gamma .$ Then
the following statements hold:

1) The algebra $\Gamma $ is Koszul with Yoneda algebra $\underset{k\geq 0}{%
\oplus }Ext_{\Gamma }^{k}(\Gamma _{0},\Gamma _{0})\cong \Lambda .$

2) If $K_{\Lambda }$ and $K_{\Gamma ^{op}}$ denote the category of Koszul $%
\Lambda $ and $\Gamma ^{op}$modules, respectively, then the ext functor
induces a duality $F$:$K_{\Lambda }$ $\rightarrow $ $K_{\Gamma ^{op}}$ such
that $F(J^{k}M\left[ k\right] )\cong \Omega ^{k}F(M)\left[ k\right] $ and $%
J^{k}F(M)\left[ k\right] \cong F(\Omega ^{k}M\left[ k\right] ).\blacksquare $
\end{theorem}

We are interested in the following examples of Koszul algebras:

1) The polynomial algebra $\Gamma=K\left[ x_{0},x_{1},...x_{r}\right] $ is a
Koszul algebra with Yoneda algebra the exterior algebra $%
\Lambda=K<x_{0},x_{1},...x_{r}>/(x_{i}^{2},$ $x_{i}x_{j}+x_{j}x_{i}).$

2) If $\Gamma$ is a Koszul $K-$algebra and $G$ a finite group of
automorphisms of $\Gamma$ such that $charK\nmid\mid G\mid$ and $\Lambda$ is
the Yoneda algebra of $\Gamma,$ then $G$ acts on $\Lambda$ and the skew
group algebra $\Gamma\ast G$ is Koszul with Yoneda algebra $\Lambda\ast G.%
\left[ M4\right] $

In the examples above the exterior algebra $\Lambda$ is selfinjective and
given a finite group of automorphisms of the $K-$ algebra $\Lambda$ such
that the $charK\nmid\mid G\mid.$ Then the skew group algebra $\Lambda\ast G$
is selfinjective if and only if $\Lambda$ is $\left[ RR\right] .$

The following theorem characterizes selfinjective Koszul algebras, it was
proved first in the connected case by P. S. Smith:

\begin{theorem}
$\left[ S1\right] ,\left[ M3\right] $Let $\Lambda $ be a finite dimensional
indecomposable Koszul algebra with Yoneda algebra $\Gamma .$ Then the
following statements are equivalent:

1) The algebra $\Lambda $ is selfinjective with radical $J$ such that $%
J^{r+1}\neq 0$ and $J^{r+2}=0.$

2) i) The graded $\Gamma ^{op}$ simple have projective dimension $r+1.$

ii) For any graded simple module $S$ the equality $Ext_{\Gamma ^{op}}^{i}($ $%
S,\Gamma )=0$ for $i\neq r+1.$

iii) The functor $Ext_{\Gamma ^{op}}^{r+1}($ $-,\Gamma )$ induces a
bijection between the graded $\Gamma $ and $\Gamma ^{op}$ simple modules.$%
\blacksquare $
\end{theorem}

\begin{definition}
We will call a graded algebra $\Gamma $ satisfying conditions i), ii), iii)
Artin Shelter regular $\left[ AS\right] $.
\end{definition}

\begin{remark}
These algebras were called generalized Auslander regular in $\left[ GMT%
\right] $, $\left[ M2\right] $, $\left[ M3\right] $
\end{remark}

Another example of Artin Shelter regular Koszul algebra is the preprojective
algebra $\Gamma=K\overset{\wedge}{Q}/I$ of a non Dynkin bipartite graph $Q,$
its Yoneda algebra is the trivial extension algebra $\Lambda=KQ\rhd D(KQ)$. $%
\left[ M1\right] $, $\left[ GMT\right] .$

We will need the following results from $\left[ M5\right] ,\left[ MM\right]
. $

\begin{theorem}
Let $\Lambda =KQ/I$ be a Koszul algebra and $M$ and $N$ two Koszul modules.
Then for any pair of integers $k$ and $l,$ with $k\geq 0$, the following two
statements are true:

i) If $Ext_{\Lambda }^{k}(M,N\left[ l\right] )_{0}\neq 0,$ then $k\geq -l.$

ii) If $k\geq -l,$ then there exists a vector space isomorphism:
\end{theorem}

\begin{center}
$Ext_{\Lambda }^{k}(M,N\left[ l\right] )_{0}\cong Ext_{\Gamma
^{op}}^{k+l}(F(N)\left[ l\right] ,F(M))_{0}.$ $\blacksquare $
\end{center}

\begin{theorem}
Let $\Gamma $ be an infinite dimensional Koszul algebra with Yoneda algebra $%
\Lambda $ and let $F:Gr\Gamma \rightarrow Gr\Lambda ^{op}$ be Koszul duality
$F(M)=\underset{k\geq 0}{\oplus }Ext_{\Gamma }^{k}(M,\Gamma _{0}).$ Then for
any pair of Koszul $\Gamma -$modules $M$ and $N,$ any integer $p$ and an
integer $n\geq 0,$ there exist a functorial isomorphism:
\end{theorem}

\begin{center}
$\underrightarrow{\lim }Ext_{\Gamma }^{n}(J^{k}M,N)_{p}\cong
\underrightarrow{\lim }Ext_{\Lambda ^{op}}^{n}(\Omega ^{k+p}F(N),\Omega
^{k}F(M))_{-p}\ \blacksquare $
\end{center}

\begin{corollary}
Under the conditions of the theorem if we assume in addition $\Lambda $
selfinjective the isomorphism becomes:
\end{corollary}

\begin{center}
$\underrightarrow{\lim }Ext_{\Gamma }^{n}(J^{k}M,N)_{p}\cong Ext_{\Lambda
^{op}}^{n}(\Omega ^{p}F(N),F(M))_{-p.}\ \blacksquare $
\end{center}

\section{Main results}

Let $\Gamma =K\left[ x_{0},x_{1},...x_{r}\right] $ be the polynomial
algebra, $\Lambda =K<x_{0},x_{1},...x_{r}>/(x_{i}^{2},$ $%
x_{i}x_{j}+x_{j}x_{i})$ the exterior algebra, $\pi :Gr\Gamma \rightarrow
QGr\Gamma $ the quotient functor. We are interested in finitely generated
graded modules $M$ such that $\pi (M)$ corresponds to a locally free sheaf
under Serre%
\'{}%
s equivalence. By the approximation proposition $\left[ M3\right] $, given a
finitely generated $\Gamma -$ module $M$ there exists an integer $k$ such
that $M_{\geq k}\left[ k\right] =N$ is Koszul. We know $\left[ P\right] ,$ $%
\left[ S2\right] $ that for two graded $\Gamma -$modules $X$ and $Y$ there
exists an integer $k$ such that $X_{\geq k}\cong Y_{\geq k}$ if and only if $%
\pi (X)\cong \pi $ $(Y).$ It follows $\pi (M)\cong \pi (N\left[ -k\right] ).$
Hence we may assume $M$ is Koszul up to shifting, this is: there exists some
integer $k$ such that $M\left[ k\right] $ is Koszul. Using Koszul duality $%
F:K_{\Lambda }\rightarrow K_{\Gamma ^{op}}$ we want to characterize the
Koszul $\Lambda -$ modules $N$ such that $\pi F(N)$ is locally free.

\begin{definition}
Let $\Gamma $ be an indecomposable noetherian Artin Shelter regular Koszul
algebra of global dimension $r+1.$ Let $\pi :Gr\Gamma \rightarrow QGr\Gamma $
be the quotient functor and M a finitely generated module. We say that $\pi $%
(M) is locally free if $\underrightarrow{\lim }$Ext$_{\Gamma }^{j}$(M$_{\geq
t}$,$\Gamma \left[ \text{s}\right] $)$_{0}$ $=0$ for $s>>0$ and $0<j\leq r.$
(See also $\left[ BGK\right] ).$
\end{definition}

Observe that if $\pi(M)$ is locally free and $k$ is an integer, then $\pi(M%
\left[ k\right] )$ is locally free:

We have: $M\left[ k\right] _{\geq t}=M_{\geq k+t}\left[ k\right] $. Set $s%
{\acute{}}%
=s-k$ and $k+t=t%
{\acute{}}%
.$

Then there exists a chain of isomorphisms:

\begin{center}
$\underrightarrow{\lim }Ext_{\Gamma }^{j}(M\left[ k\right] _{\geq t},\Gamma %
\left[ s\right] )_{0}\cong \underrightarrow{\lim }Ext_{\Gamma }^{j}(M_{\geq
k+t}\left[ k\right] ,\Gamma \left[ s\right] )_{0}\cong $

$\underrightarrow{\lim }Ext_{\Gamma }^{j}(M_{\geq k+t},\Gamma \left[ s-k%
\right] )_{0}\cong \underrightarrow{\lim }Ext_{\Gamma }^{j}(M_{\geq t%
{\acute{}}%
},\Gamma \left[ s%
{\acute{}}%
\right] )_{0}=0$
\end{center}

for $s^{\prime }>>0$ and $0<j\leq r.$

We remarked above that for any finitely generated module $M$ there exists
some integer $k$ such that $M_{\geq k}\left[ k\right] $ is Koszul. We are
interested in Koszul $\Gamma^{op}-$ modules $M$ such that $\pi(F(M))$ is
locally free.

\begin{theorem}
Let $\Lambda $ be an indecomposable selfinjective Koszul algebra with $%
J^{r+1}\neq 0$ and $J^{r+2}=0.$ Assume the Yoneda algebra $\Gamma $ is
noetherian. Let $F:K_{\Lambda }\rightarrow K_{\Gamma ^{op}}$ be Koszul
duality. Then a Koszul $\Lambda -$ module $N$ has the property that $\pi
(F(N))$ is locally free if and only if $Ext_{\Lambda }^{j}(\Omega
^{s}\Lambda _{0}\left[ s\right] ,N)_{0}=0$ for $0<j\leq r$ and $s>>0.$
\end{theorem}

\textit{Proof. }By Theorem 7, we have natural isomorphisms:

\begin{center}
$\underrightarrow{\lim }Ext_{\Gamma ^{op}}^{j}(J^{k}F(N),\Gamma \left[ s%
\right] )_{0}\cong \underrightarrow{\lim }Ext_{\Gamma
^{op}}^{j}(J^{k}F(N),\Gamma )_{s}\cong $

$\underrightarrow{\lim }Ext_{\Lambda }^{j}(\Omega ^{k+s}F^{-1}(\Gamma
),\Omega ^{k}N)_{-s}\cong \underrightarrow{\lim }Ext_{\Lambda }^{j}(\Omega
^{k+s}\Lambda _{0}\left[ s\right] ,\Omega ^{k}N))_{0}$\\[0pt]
$\cong Ext_{\Lambda }^{j}(\Omega ^{s}\Lambda _{0}\left[ s\right] ,N)_{0.}$
\end{center}

It follows $\pi (F(N))$ is locally free if and only if $Ext_{\Lambda
}^{j}(\Omega ^{s}\Lambda _{0}\left[ s\right] ,N)_{0}=0$ for $0<j\leq r$ and $%
s>>0.$ $\blacksquare $

\begin{corollary}
Let $N$ be a Koszul $\Lambda -$module. Then for any positive integer $k$ the
module $\Omega ^{k}N\left[ k\right] $ is such that $\pi (F(\Omega ^{k}N\left[
k\right] ))$ is locally free if and only if $\pi (F(N))$ is locally free.
\end{corollary}

\textit{Proof. }$Ext_{\Lambda }^{j}(\Omega ^{s}\Lambda _{0}\left[ s\right]
,\Omega ^{k}N\left[ k\right] )_{0}\cong Ext_{\Lambda }^{j}(\Omega
^{s-k}\Lambda _{0}\left[ s-k\right] ,N)_{0}=0$ for $s>>0.\blacksquare $

\begin{lemma}
Let $\Lambda =KQ/I$ be a graded quiver algebra, let $M$ and $N$ be graded
modules generated in degree $l$ and $k$, respectively, with $l>k.$ Then
there is a natural isomorphism:
\end{lemma}

\begin{center}
$Ext_{\Lambda }^{1}(M,N)_{0}\cong Ext_{\Lambda }^{1}(M,J^{l-k}N)_{0}.$
\end{center}

\textit{Proof. }The exact sequence: $0\rightarrow J^{l-k}N\rightarrow
N\rightarrow N/J^{l-k}N\rightarrow 0$ induces an exact sequence:

\begin{center}
$Hom_{\Lambda }(M,N/J^{l-k}N)_{0}\rightarrow Ext_{\Lambda
}^{1}(M,J^{l-k}N)_{0}\overset{\delta }{\rightarrow E}xt_{\Lambda
}^{1}(M,N)_{0}.$
\end{center}

Where $Hom_{\Lambda }(M,N/J^{l-k}N)_{0}=0.$

We must prove $\delta $ is an epimorphism.

Let $x\in Ext_{\Lambda }^{1}(M,N)_{0}.$ Then the exact sequence $x:$ $%
0\rightarrow N\rightarrow E\rightarrow M\rightarrow 0$ induces an exact
commutative diagram:

\begin{center}
$
\begin{array}{ccccccc}
0\rightarrow & N_{\geq l} & \rightarrow & E_{\geq l} & \rightarrow & M &
\rightarrow 0 \\
& \downarrow &  & \downarrow &  & \downarrow 1 &  \\
0\rightarrow & N & \rightarrow & E & \rightarrow & M & \rightarrow 0%
\end{array}
$
\end{center}

With $N_{\geq l}\cong J^{l-k}N$ and we have proved $\delta $ is an
epimorphism. $\blacksquare $

\begin{proposition}
Let $\Lambda $ be an indecomposable selfinjective Koszul algebra with $%
J^{r+1}\neq 0$ and $J^{r+2}=0.$ Assume the Yoneda algebra $\Gamma $ is
noetherian and let $F:K_{\Lambda }\rightarrow K_{\Gamma ^{op}}$ be the
Koszul duality. Let $0\rightarrow X\rightarrow Y\rightarrow Z\rightarrow 0$
be an exact sequence of Koszul $\Lambda -$ modules and degree zero maps.
Then the following statements hold:

1) If $\pi (F(X))$ and $\pi (F(Z))$ are locally free, then $\pi (F(Y))$ is
locally free.

2) If $\pi (F(X))$ and $\pi (F(Y))$ are locally free, then $\pi (F(Z))$ is
locally free.
\end{proposition}

\textit{Proof. }1) The exact sequence $0\rightarrow X\rightarrow
Y\rightarrow Z\rightarrow 0$ of Koszul $\Lambda -$ modules induces an exact
sequence:

\begin{center}
$Ext_{\Lambda }^{i}(\Omega ^{s}\Lambda _{0}\left[ s\right]
,X)_{0}\rightarrow Ext_{\Lambda }^{i}(\Omega ^{s}\Lambda _{0}\left[ s\right]
,Y)_{0}\rightarrow Ext_{\Lambda }^{i}(\Omega ^{s}\Lambda _{0}\left[ s\right]
,Z)_{0},$
\end{center}

by hypothesis, $Ext_{\Lambda }^{i}(\Omega ^{s}\Lambda _{0}\left[ s\right]
,X)_{0}=Ext_{\Lambda }^{i}(\Omega ^{s}\Lambda _{0}\left[ s\right] ,Z)_{0}=0$
for $0<i\leq r$\linebreak\ and $s>>0.$

Hence; $Ext_{\Lambda }^{i}(\Omega ^{s}\Lambda _{0}\left[ s\right] ,Y)_{0}=0$
for $0<i\leq r$ and $s>>0.$

2) Let $k$ be an integer $k\geq r+1$ and $L$ a torsion free $\Gamma -$
module.

It was proved in $\left[ M2\right] ,$ $Ext_{\Gamma ^{op}}^{k}(L,\Gamma )=0,$
hence; $Ext_{\Gamma ^{op}}^{k}(L,\Gamma )_{j}=0$ for all $j$ and $k\geq r+1.$

If $F(X)$ is torsion free, then using lemma 1 and theorem 6, we have a chain
of natural isomorphisms:

\begin{center}
$Ext_{\Lambda }^{k}(\Omega ^{s}\Lambda _{0}\left[ s\right] ,X)_{0}\cong
Ext_{\Lambda }^{1}(\Omega ^{k+s-1}\Lambda _{0}\left[ s\right] ,X)_{0}\cong
Ext_{\Lambda }^{1}(\Omega ^{k+s-1}\Lambda _{0}\left[ s\right] ,J^{k-1}X)_{0}$

$\cong Ext_{\Gamma op}^{1}(\Omega ^{k-1}F(X)\left[ -s\right]
,J^{k+s-1}F(\Lambda _{0}))_{0}\cong Ext_{\Gamma op}^{1}(\Omega ^{k-1}F(X)%
\left[ -s\right] ,F(\Lambda _{0}))_{0}$

$\cong Ext_{\Gamma ^{op}}^{k}(F(X),F(\Lambda _{0}))_{s}=Ext_{\Gamma
^{op}}^{k}(F(X),\Gamma )_{s}=0$ for $k\geq r+1.$
\end{center}

Assume $F(X)$ has torsion and let $S_{1}$ be a simple in the socle, then
there is an integer $k$ such that $J^{k}F(X)\cong F(\Omega ^{k}X)\cong
S_{1}\oplus W_{1}.$ Applying $F^{-1},$ there is an isomorphism: $\Omega
^{k}X\cong P_{1}\oplus X_{1},$ but since $P_{1}$ is projective injective $%
k=0 $ and $S_{1}$ appears at the top of $F(X).$ We have proved $F(X)\cong
S\oplus W$ with $S$ semisimple and $W$ torsion free, hence $X\cong P\oplus
X^{\prime }$ with $P$ projective and $X^{\prime }$ is such that $\pi
(F(X^{\prime }))$ is locally free.

We have an exact sequence:

\begin{center}
$Ext_{\Lambda }^{i}(\Omega ^{s}\Lambda _{0}\left[ s\right]
,Y)_{0}\rightarrow Ext_{\Lambda }^{i}(\Omega ^{s}\Lambda _{0}\left[ s\right]
,Z)_{0}\rightarrow Ext_{\Lambda }^{i+1}(\Omega ^{s}\Lambda _{0}\left[ s%
\right] ,X)_{0}$.
\end{center}

By hypothesis and the above remark, for $0<i\leq r$ and $s>>0$ we have:

\begin{center}
$Ext_{\Lambda }^{i}(\Omega ^{s}\Lambda _{0}\left[ s\right]
,Y)_{0}=0=Ext_{\Lambda }^{i+1}(\Omega ^{s}\Lambda _{0}\left[ s\right]
,X^{\prime })_{0}$ $.$
\end{center}

Hence; $Ext_{\Lambda }^{i}(\Omega ^{s}\Lambda _{0}\left[ s\right] ,Z)_{0}=0$
for $0<i\leq r$ and $s>>0.$ $\blacksquare $

\begin{proposition}
Let $\Lambda $ be an indecomposable selfinjective Koszul algebra with $%
J^{r+1}\neq 0$ and $J^{r+2}=0.$ Assume the Yoneda algebra $\Gamma $ is
noetherian and let $F:K_{\Lambda }\rightarrow K_{\Gamma ^{op}}$ be the
Koszul duality. Then the following statement holds:

If $X$ is a Koszul $\Lambda -$ module such that $\pi (F(X))$ is locally
free, then $JX\left[ 1\right] $ is such that $\pi (F(JX\left[ 1\right] ))$
is locally free.
\end{proposition}

\textit{Proof. }Let $X$ be a Koszul $\Lambda -$ module such that $\pi (F(X))$
is locally free, let $0\rightarrow \Omega (X)$ $\rightarrow P\rightarrow
X\rightarrow 0$ be exact with $P$ the projective cover of $X.$ It induces an
exact sequence: $0\rightarrow \Omega (X)\left[ 1\right] $ $\rightarrow JP%
\left[ 1\right] \rightarrow JX\left[ 1\right] \rightarrow 0$. By Corollary
2, $\pi (F(\Omega (X)\left[ 1\right] )$ is locally free.

If $P$ is projective generated in degree zero and $M$ is a Koszul module,
then from the exact sequence: $0\rightarrow JP\left[ 1\right] \rightarrow P%
\left[ 1\right] \rightarrow S\left[ 1\right] \rightarrow 0$ we obtain an
exact sequence:

\begin{center}
$0\rightarrow Hom_{\Lambda }(M,JP\left[ 1\right] )_{0}\rightarrow
Hom_{\Lambda }(M,P\left[ 1\right] )_{0}\rightarrow $

$\rightarrow Hom_{\Lambda }(M,S\left[ 1\right] )_{0}\rightarrow Ext_{\Lambda
}^{1}(M,JP\left[ 1\right] )_{0}\rightarrow 0$
\end{center}

and isomorphisms:

\begin{center}
$Ext_{\Lambda }^{j-1}(M,S\left[ 1\right] )_{0}\cong Ext_{\Lambda }^{j}(M,JP%
\left[ 1\right] )_{0}.$
\end{center}

Since $M$ is generated in degree zero and $\Omega ^{j-1}M$ in degree $j-1,$
it follows $Hom_{\Lambda }(M,S\left[ 1\right] )_{0}=0$ and $Ext_{\Lambda
}^{j-1}(M,S\left[ 1\right] )_{0}=0.$ Therefore $Ext_{\Lambda }^{j}(M,JP\left[
1\right] )_{0}=0$ for $j>0.$

Letting $M$ be equal to $\Omega ^{s}\Lambda _{0}\left[ s\right] $ we get $%
Ext_{\Lambda }^{i}(\Omega ^{s}\Lambda _{0}\left[ s\right] ,JP\left[ 1\right]
)_{0}=0$ for all $s>0$ and all $i>0.$

We have proved $\pi (F(JP\left[ 1\right] ))$ is locally free. It follows by
2) of previous proposition $\pi (F(JX\left[ 1\right] ))$ is locally free. $%
\blacksquare \smallskip $

We will consider now the case $r=1.$

\begin{theorem}
Let $\Gamma $ be an Artin Shelter regular Koszul algebra of global dimension
2 and assume $\Gamma $ is noetherian indecomposable with Yoneda algebra $%
\Lambda ,$ let $F:K_{\Lambda }\rightarrow K_{\Gamma ^{op}}$ be the Koszul
duality. A Koszul $\Lambda -$ module $M$ is such that $\pi (F(M))$ is
locally free if and only if $M\cong \Omega ^{s}S\left[ s\right] $ for a
simple module $S$ and some integer $s\geq 0.$
\end{theorem}

\textit{Proof. }A Koszul $\Lambda -\mathit{mod}$ule $M$ is such that $\pi
(F(M))$ is locally free if and only if $Ext_{\Lambda }^{1}(\Omega
^{s}\Lambda _{0}\left[ s\right] ,M)_{0}=0$ for $s>>0.$ Equivalently, if and
only if the maps module projectives \underline{$Hom$}$_{\Lambda }(\Omega
^{s+1}\Lambda _{0}\left[ s\right] ,M)_{0}=0$ for $s>>0$ .

Assume $M=\Omega ^{t}S\left[ t\right] ,$ with $S$ a simple $\Lambda -$
module and let $s>t.$

Then \underline{$Hom$}$_{\Lambda }(\Omega ^{s+1}\Lambda _{0}\left[ s\right] $%
,$\Omega ^{t}$S$\left[ t\right] )_{0}=$\underline{$Hom$}$_{\Lambda }(\Omega
^{s-t+1}\Lambda _{0}\left[ s\text{-}t\right] $,$S)_{0},$ but $\Omega
^{s-t+1} $S$\left[ s\text{-}t\right] $ is generated in degree $1$, hence;
\underline{$Hom$}$_{\Lambda }(\Omega ^{s-t+1}\Lambda _{0}\left[ s-t\right]
,S)_{0}=0.$

Conversely, assume $M$ is not isomorphic to $\Omega ^{s}S\left[ s\right] $
for any simple $S$ and any $s.$

There is an isomorphism: \underline{$Hom$}$_{\Lambda }(\Omega ^{s+1}\Lambda
_{0}\left[ s\right] ,M)_{0}=$\underline{$Hom$}$_{\Lambda }(\Lambda
_{0},\Omega ^{-s-1}M\left[ -s\right] )_{0}.$

If the Loewy length of $M$ is two, then $socM$ is generated in degree one.
We have an exact sequence: $0\rightarrow M\rightarrow P\overset{p}{%
\rightarrow }\Omega ^{-1}M\rightarrow 0$ with $P$ generated in degree $-1.$
Since $\Omega ^{-1}M$ is not projective, $p(socP)=0$ and Loewy length of $%
\Omega ^{-1}M$ $\leq 2.$

If $\Omega ^{-1}M$ is not simple, then $soc(\Omega ^{-1}M$ )$\left[ -1\right]
$ is generated in degree $1.$ It follows by induction $soc\Omega ^{-s}M\left[
-s\right] $ is generated in degree $1,$ hence; $\Omega ^{-s-1}M\left[ -s%
\right] $ has socle in degree zero and \underline{$Hom$}$_{\Lambda }(\Lambda
_{0},\Omega ^{-s-1}M\left[ -s\right] )_{0}\neq 0$ for all $s.\blacksquare
\smallskip $

In the general case we have the following remarks:

\begin{proposition}
Let $\Gamma $ be an indecomposable Artin Shelter regular noetherian Koszul
algebra of global dimension $r+1$ with Yoneda algebra $\Lambda .$ Let $%
F:K_{\Lambda }\rightarrow K_{\Gamma ^{op}}$ be the Koszul duality. Let $N$
be a Koszul $\Lambda -$ module such that $\pi (F(N))$ is locally free. Then
there exists some integer $s$ such that $\Omega ^{-s}N\left[ -s\right] $ is
not Koszul.
\end{proposition}

\textit{Proof. }Assume for all $s\geq 0$ the module $\Omega ^{-s}N\left[ -s%
\right] $ is Koszul. The fact $\pi (F(N))$ is locally free implies for $s>>0$
and $0<j\leq r$ there are isomorphisms:

\begin{center}
\underline{$Hom$}$_{\Lambda }(\Omega ^{s+j}\Lambda _{0}\left[ s\right]
,N)_{0}\cong $\underline{$Hom$}$_{\Lambda }(\Omega ^{j}\Lambda _{0},\Omega
^{-s}N\left[ -s\right] )_{0}$

$\cong $\underline{$Hom$}$_{\Lambda }(\Lambda _{0},\Omega ^{-s-j}N\left[ -s%
\right] )_{0}=0$
\end{center}

By hypothesis the module $X=\Omega ^{-s-j}N\left[ -s-j\right] $ is Koszul.
Consider the exact sequence: $0\rightarrow X\rightarrow P\rightarrow \Omega
^{-1}X\rightarrow 0$ with $P$ the injective envelope of X. Since $\Omega
^{-1}X\left[ -1\right] $ is Koszul then $SocX$ is generated in degree $r.$

If $j=r,$ then $X\left[ j\right] $ has socle generated in degree zero and
\underline{$Hom$}$_{\Lambda }(\Lambda _{0},X\left[ j\right] )_{0}=$

\underline{$Hom$}$_{\Lambda }(\Lambda _{0},\Omega ^{-s-j}N\left[ -s\right]
)_{0}\neq 0$ is a contradiction.$\blacksquare \smallskip $

We have the following characterization of the Koszul $\Lambda-$ modules $M$
with $\pi(F(M))$ locally free:

\begin{theorem}
Let $\Gamma $ be an indecomposable noetherian Artin Shelter regular Koszul
algebra of global dimension $r+1$ with Yoneda algebra $\Lambda ,$ let $%
F:K_{\Lambda }\rightarrow K_{\Gamma ^{op}}$ be the Koszul duality. Then an
indecomposable Koszul module $M$ is such that $\pi (F(M))$ is locally free
if and only if there exists an integer $s>0$ such that $\Omega ^{-s}M\left[
-s\right] _{\geq 1}=0.$
\end{theorem}

\textit{Proof. }Assume $\pi (F(M))$ is locally free and suppose for all $s>0$
the module $\Omega ^{-s}M\left[ -s\right] _{\geq 1}\neq 0.$

We know by $\left[ MZ\right] ,$ there exists an integer $s$ such that $%
\Omega ^{-s}M\left[ -s\right] $ is weakly co-Koszul.

Let $X=\Omega ^{-s}M\left[ -s\right] $ and let $0\rightarrow X\rightarrow
Q\rightarrow \Omega ^{-1}X\rightarrow 0$ be exact with $Q$ the injective
envelope of $X.$ Since $X$ is weakly co-Koszul the epimorphism: $%
Q/socQ\rightarrow \Omega ^{-1}X\rightarrow 0$ induces an epimorphism: $%
soc^{2}Q/socQ\rightarrow soc\Omega ^{-1}X\rightarrow 0$ this means that if $%
X $ has a cogenerator in degree $k,$ then $\Omega ^{-1}X$ has a cogenerator
in degree $k-1$ and $\Omega ^{-1}X\left[ -1\right] $ has a cogenerator in
degree $k.$

The module $M$ is generated in degree zero, then it has cogenerators in
degree $k$ with $k\leq r$ and $\Omega ^{-1}M\left[ -1\right] $ has
cogenerators in degree $\leq r,$ it follows by induction $X$ has
cogenerators in degree $\leq r,$ by hypothesis $\Omega ^{-s}M\left[ -s\right]
_{\geq 1}\neq 0,$ then $\Omega ^{-s}M\left[ -s\right] $ has a cogenerator in
degree $k$ with $1\leq k\leq r$ $.$ It follows for all $t\geq 0$ the module $%
\Omega ^{-t-s}M\left[ -t-s\right] $ has cogenerator in lowest degree $1\leq
k\leq r.$ Then $\Omega ^{-k}\Omega ^{-(s+t)}M\left[ -(s+t)\right] $ has a
cogenerator in degree zero and it follows:

\underline{$Hom$}$_{\Lambda }(\Lambda _{0},\Omega ^{-k}\Omega ^{-(s+t)}M%
\left[ -(s+t)\right] )_{0}\neq 0,$ hence; $Ext_{\Lambda }^{k}(\Omega ^{\ell
}\Lambda _{0}\left[ l\right] ,M)_{0}\neq 0$ for all $\ell \geq s$ a
contradiction.

Assume now there exists an $s$ such that $\Omega ^{-s}M\left[ -s\right]
_{\geq 1}=0$ and set $X$=$\Omega ^{-s}M\left[ -s\right] $.

As above, if $X$ has a cogenerator in lowest degree $k,$ then $\Omega ^{-1}X%
\left[ -1\right] $ has a cogenerator in lowest degree $\leq k.$ It follows
for all $t\geq 0$ the module $\Omega ^{-(s+t)}$M$\left[ \text{-s-t}\right]
_{\geq 1}=0$.

Hence; for all $0<j\leq r$ the module $\Omega ^{-j}\Omega ^{-(s+t)}M\left[
-(s+t)\right] $ has cogenerators in degree less than zero. It follows: $%
Hom_{\Lambda }(\Lambda _{0},\Omega ^{-j}\Omega ^{-(s+t)}M\left[ -(s+t)\right]
)_{0}=0$ for all $t\geq 0$ and $0<j\leq r.$

We have proved $\pi (F(M))$ is locally free.$\blacksquare $

\begin{corollary}
Let $\Lambda $ and $\Gamma $ be as in the theorem. Let $M$ be a Koszul and
quasi co-Koszul $\Lambda -$ module. Then $\pi (F(M))$ is locally free if and
only if $M$ is simple.
\end{corollary}

\textit{Proof. }It is clear that if $M$ is simple, then $F(M)$ is projective
and $\pi (F(M))$ locally free.

Assume $M$ is not simple. Then $M$ has cogenerators in degree $k$ with $%
r\geq k>0,$ hence; $M_{\geq 1}\neq 0.$

As in the proof of the theorem, $\Omega ^{-s}M\left[ -s\right] $ has
cogenerators in degree $k>0.$ Therefore: $Hom_{\Lambda }(\Lambda _{0},\Omega
^{-k}\Omega ^{-s}M\left[ -s\right] )_{0}\neq 0$ for all $s>0.$

It follows $\pi (F(M))$ is not locally free.$\blacksquare \smallskip $

We can improve now the characterization of the $\Lambda-$modules
corresponding to locally free sheaves.

\begin{theorem}
Let $\Gamma $ be an indecomposable noetherian Artin Schelter regular Koszul
algebra of global dimension $r+1$ with Yoneda algebra $\Lambda ,$ let $%
F:K_{\Lambda }\rightarrow K_{\Gamma ^{op}}$ be Koszul duality. Then an
indecomposable Koszul module $M$ is such that $\pi F(M)$ is locally free if
and only if there exists a non negative integer $t$ such that $\Omega ^{-t}M%
\left[ -t\right] $ is co-Koszul.
\end{theorem}

\textit{Proof. }Assume there exists an integer $t\geq 0$ such that $\Omega
^{-t}M\left[ -t\right] $ is co-Koszul. \linebreak By hypothesis, $\Omega
^{-t}M\left[ -t\right] $ is cogenerated in degree zero, this
means\linebreak\ $\Omega ^{-t}M\left[ -t\right] $ $_{\geq 1}=0$. Applying
Theorem 10, it follows $\pi F(M)$ is locally free.

We assume now $\pi F(M)\;$is locally free, by Theorem 10, there exists an
integer $t\geq 0$ such that $\Omega ^{-t}M\left[ -t\right] $ $_{\geq 1}=0.$
Let $S$ be a simple in the socle of $\Omega ^{-t}M\left[ -t\right] $. By
hypothesis, $S$ is generated in degree $k\leq 0.$ Suppose $k<0,$ then $%
\Omega ^{t}S$ is generated in degree $k+t$ and $M\left[ -t\right] $ is
generated in degree $t$ with $k+t<t.$ It follows $Hom_{\Lambda }(\Omega
^{t}S,M\left[ -t\right] )_{0}=0.$

We have isomorphisms: $Hom_{\Lambda }(S,\Omega ^{-t}M\left[ -t\right] )_{0}=%
\underline{Hom}_{\Lambda }(S,\Omega ^{-t}M\left[ -t\right] )_{0}\cong $%
\linebreak $\underline{Hom}_{\Lambda }(\Omega ^{t}S,M\left[ -t\right] )_{0}$
and $Hom_{\Lambda }(S,\Omega ^{-t}M\left[ -t\right] )_{0}\neq 0$. Since the
natural inclusion: $j:S\rightarrow \Omega ^{-t}M\left[ -t\right] $ is not
zero, we have reached a contradiction, proving \linebreak $soc\Omega ^{-t}M%
\left[ -t\right] $ is generated in degree zero.

Let $s>t.$ It is clear $\Omega ^{-s}M\left[ -s\right] $ $_{\geq 1}=0,$
hence; $soc\Omega ^{-s}M\left[ -s\right] $ is generated in degree zero. It
follows $\Omega ^{-t}M\left[ -t\right] $ is co-Koszul.$\blacksquare $

\begin{theorem}
Let $\Gamma $and $\Lambda $ be Koszul algebras satisfying the conditions of
Theorem 11 and $F:K_{\Lambda }\rightarrow K_{\Gamma ^{op}}$ be Koszul
duality. Let $M$ be a Koszul $\Lambda -$module such that $\pi F(M)$ is
locally free. Then there exists a semisimple module $S$ finitely generated
in degree zero, a non negative integer $t$ and an epimorphism: $\Omega ^{t}S%
\left[ t\right] \rightarrow M\rightarrow 0.$
\end{theorem}

\textit{Proof. }Let $t$ be a non negative integer such that $\Omega ^{-t}M%
\left[ -t\right] $ is co-Koszul and let $S$ be the socle of $\Omega ^{-t}M%
\left[ -t\right] $ $,$the inclusion: $j:S\rightarrow \Omega ^{-t}M\left[ -t%
\right] $ induces a non zero map: $\Omega ^{t}j:\Omega ^{t}S\rightarrow M%
\left[ -t\right] .$We want to prove that $\Omega ^{t}j$ is an epimorphism.

Assume $\Omega ^{t}j$ is not an epimorphism, then there exists a simple $%
S^{\prime }$ generated in degree zero and a non zero map $h:M\left[ -t\right]
\rightarrow S^{\prime }\left[ -t\right] $ such that $h$ $\Omega ^{t}j=0.$
Therefore: $(\Omega ^{-t}h)j=0$ and $Ker\Omega ^{-t}h\supseteq soc\Omega
^{-t}M\left[ -t\right] .$

Let $T$ be a simple module contained in $soc(Im\Omega ^{-t}h),$taking the
pull back of the inclusion $i$ we obtain a commutative exact
diagram:\medskip \newline
$
\begin{array}{ccccccc}
\text{0}\rightarrow & \text{Ker}\Omega ^{-t}\text{h/soc}\Omega ^{-t}\text{M}%
\left[ \text{-t}\right] & \rightarrow & \text{W} & \rightarrow & \text{T} &
\rightarrow \text{0} \\
& \downarrow \text{1} &  & \downarrow &  & \downarrow \text{i} &  \\
\text{0}\rightarrow & \text{Ker}\Omega ^{-t}\text{h/soc}\Omega ^{-t}\text{M}%
\left[ \text{-t}\right] & \rightarrow & \Omega ^{-t}\text{M}\left[ -t\right]
\text{/soc}\Omega ^{-t}\text{M}\left[ \text{-t}\right] & \rightarrow & \text{%
Im}\Omega ^{-t}\text{h} & \rightarrow \text{0}%
\end{array}
\smallskip \smallskip $

Since $\Omega ^{-t}M\left[ -t\right] $ is cogenerated in degree zero $%
W_{k}\neq 0$ implies $k<0$. Therefore: $W$ has all its generators in
negative degrees. It follows $T$ is generated in negative degree, but this
is a contradiction since $soc(Im\Omega ^{-t}h)\subseteq $ $soc\Omega
^{-t}S^{\prime }\left[ -t\right] $ and $soc\Omega ^{-t}S^{\prime }\left[ -t%
\right] $ is generated in degree zero.$\blacksquare $

\section{The Auslander Reiten quiver of K$_{\Lambda }.$}

In this section, we study the Auslander Reiten quiver of the category $%
K_{\Lambda }$ of Koszul modules over a selfinjective Koszul algebra $\Lambda
.$We recall from $\left[ GMRSZ\right] $, that if $M$ is an indecomposable
Koszul module, then there exists an Auslander sequence: $0\rightarrow \sigma
(M)\rightarrow E\rightarrow M\rightarrow 0$ in $K_{\Lambda }$ ending at $M$.
This sequence is constructed as follows: we first look at the Auslander
Reiten sequence ending at $M$ in $gr_{\Lambda }:0\rightarrow \tau
(M)\rightarrow F\rightarrow M\rightarrow 0$ and then, if we define $\sigma
(M)=(\tau M)_{\geq 0}$ and $E=F_{\geq 0}$, it turns out that $\sigma (M)$ is
an indecomposable Koszul module, and we also get that the induced sequence
is an Auslander-Reiten sequence in $K_{\Lambda .}$ Thus $K_{\Lambda }$ has
left A-R sequences. Also note that $\sigma (M)$ has Loewy length exactly
two, hence it is also a $\Lambda /J^{2}-$module. We also mention that if $%
0\rightarrow A\rightarrow B\rightarrow C\rightarrow 0$ is an A-R sequence in
$K_{\Lambda }$, then it is also an A-R sequence in the category of graded
modules finitely generated in degree zero, $gr_{0\Lambda }$. This will be
used through the section.

\begin{lemma}
Let $0\rightarrow A\rightarrow B\rightarrow C\rightarrow 0$ be a non split
exact sequence in $K_{\Lambda }.$Then the induced sequence $0\rightarrow
\sigma (A)\rightarrow \sigma (B)\rightarrow \sigma (C)\rightarrow 0$ is also
exact in $K_{\Lambda }$.
\end{lemma}

\textit{Proof. }Since for each $i\geq 0$ we have $J^{i}A=J^{i}B\cap A$, we
know from $\left[ MZ\right] $, that the exactness of $0\rightarrow
A\rightarrow B\rightarrow C\rightarrow 0$ implies that we have an induced
exact sequence: $0\rightarrow \tau (A)\rightarrow \tau (B)\rightarrow \tau
(C)\rightarrow 0$. The lemma follows immediately from the definition of $%
\sigma .\blacksquare \medskip $

The main tool that will be used very often, is the following result from $%
\left[ GMRSZ\right] :$

\begin{lemma}
Let $\mathcal{L}_{\Lambda }$ denote the full subcategory of $gr_{0\Lambda }$
consisting of the modules $M$ having a linear presentation, that is a
projective presentation of the form: $P_{1}\rightarrow P_{0}\rightarrow
M\rightarrow 0$ where $P_{1}$ is generated in degree $1$ and $P_{0}$ is
generated in degree zero$.$ Then $\Lambda /J^{2}\otimes _{\Lambda }-:%
\mathcal{L}_{\Lambda }\rightarrow gr_{0\Lambda /J^{2}}$ is an equivalence of
categories. $\blacksquare $\newline
\end{lemma}

The following results are true in a more general setting for certain
subcategories of the graded module category of a finite dimensional algebra
having left Auslander-Reiten sequences. They are part of the finite
dimensional algebra folklore and we include proofs only for our more
restricted needs.

\begin{proposition}
Let $M$ be an indecomposable non projective Koszul $\Lambda -$module and let
$0\rightarrow \sigma (M)\overset{f}{\rightarrow }E\overset{g}{\rightarrow }%
M\rightarrow 0$ be the A-R sequence in $K_{\Lambda }$ ending at $M.$ Let $f=%
\left[ f_{1},f_{2},...f_{k}\right] $ and $g=$ $\left[ g_{1},g_{2},...g_{k}%
\right] ^{T},$where $E=E_{1}\oplus E_{2}\oplus ...E_{k},$and for each $1\leq
i\leq k$, $E_{i}$ is indecomposable in $K_{\Lambda }.$ Then, for each $1\leq
i\leq k$, the maps $f_{i}$ and $g_{i}$ are irreducible in $K_{\Lambda }.$
\end{proposition}

\textit{Proof. }First note that, by construction, the sequence\allowbreak\ $%
0\rightarrow \sigma (M)\overset{f}{\rightarrow }E\overset{g}{\rightarrow }%
M\rightarrow 0$ is also an A-R sequence in $gr_{0\Lambda }$. We show, more
generally, that the maps are in fact irreducible in $gr_{0\Lambda }.$

Assume now we have a factorization of $g_{i}$ in $gr_{0\Lambda }$:$
\begin{array}{ccc}
\text{E}_{i} & \overset{g_{i}}{\rightarrow } & \text{M} \\
\underset{\text{j}}{\searrow } &  & \underset{\text{h}}{\nearrow } \\
& \text{X} &
\end{array}
$ and that h is a not splittable epimorphism.

Then we have an induced factorization

\begin{center}
$
\begin{array}{cccc}
\text{E}_{i}\oplus \text{E}_{i}^{\prime } & \overset{\left[
g_{i}.g_{i}^{\prime }\right] }{\longrightarrow } &  & \text{M} \\
\underset{\left[
\begin{array}{cc}
\text{j} & 0 \\
0 & \text{1}%
\end{array}
\right] }{\searrow } &  & \underset{\left[ h,g_{i}^{\prime }\right] }{%
\nearrow } &  \\
& \text{X}\oplus \text{E}_{i}^{\prime } &  &
\end{array} $ \end{center} 
where $E_{i}^{\prime }=\underset{i\neq s}{\oplus }E_{s}$. Lets us show that $%
\left[ h,g_{i}^{\prime }\right] $ is not a splittable epimorphism. If it is,
there exists a morphism$\left[
\begin{array}{c}
s \\
t%
\end{array}
\right] :M\rightarrow X\oplus E_{i}^{\prime }$, such that $\left[
h,g_{i}^{\prime }\right] \left[
\begin{array}{c}
s \\
t%
\end{array}
\right] =1_{M}$, and we have $hs+g_{i}^{\prime }t=1_{M}.$ The composition $%
g_{i}^{\prime }t$ is an endomorphism of M, and since $g_{i}^{\prime }$ is
not a splittable epimorphism, the image is in the radical of $EndM$, hence $%
g_{i}^{\prime }t$ is nilpotent. Therefore $hs=1-g_{i}^{\prime }t$ is
invertible in $EndM$, so $h$ is a splittable epimorphism contradicting our
assumption on $h.$ This means that we can lift $\left[ h,g_{i}^{\prime }%
\right] $ to $E_{i}\oplus E_{i}^{\prime }$, and we obtain the following
commutative diagram with exact rows:

\begin{center}
$
\begin{array}{ccccccc}
0\rightarrow & \sigma (M) & \rightarrow & E_{i}\oplus E_{i}^{\prime } &
\overset{\left[ g_{i},g_{i}^{\prime }\right] }{\rightarrow } & M &
\rightarrow 0 \\
& \downarrow &  & \downarrow \left[
\begin{array}{cc}
j & 0 \\
0 & 1%
\end{array}
\right] &  & \parallel &  \\
0\rightarrow & K & \rightarrow & X\oplus E_{i}^{\prime } & \overset{\left[
h,g_{i}^{\prime }\right] }{\rightarrow } & M & \rightarrow 0 \\
& \downarrow &  & \downarrow \left[
\begin{array}{cc}
a & b \\
c & d%
\end{array}
\right] &  & \parallel &  \\
0\rightarrow & \sigma (M) & \rightarrow & E_{i}\oplus E_{i}^{\prime } &
\overset{\left[ g_{i},g_{i}^{\prime }\right] }{\rightarrow } & M &
\rightarrow 0%
\end{array}
$
\end{center}

In the composite diagram

\begin{center}
$
\begin{array}{ccccccc}
0\rightarrow & \sigma (M) & \rightarrow & E_{i}\oplus E_{i}^{\prime } &
\overset{\left[ g_{i},g_{i}^{\prime }\right] }{\rightarrow } & M &
\rightarrow 0 \\
& \downarrow \alpha &  & \downarrow \left[
\begin{array}{cc}
aj & b \\
cj & d%
\end{array}
\right] &  & \parallel &  \\
0\rightarrow & \sigma (M) & \rightarrow & E_{i}\oplus E_{i}^{\prime } &
\overset{\left[ g_{i},g_{i}^{\prime }\right] }{\rightarrow } & M &
\rightarrow 0%
\end{array}
$
\end{center}

we observe that $\alpha $ cannot be nilpotent, otherwise a simple argument
would show that $1_{M}$ factors through $E_{i}\oplus E_{i}^{\prime }$,
contradicting our assumptions. Since $\sigma (M)$ is indecomposable, $\alpha
$ must be invertible. Thus $\left[
\begin{array}{cc}
aj & b \\
cj & d%
\end{array}
\right] =\left[
\begin{array}{cc}
a & b \\
c & d%
\end{array}
\right] \left[
\begin{array}{cc}
j & 0 \\
0 & 1%
\end{array}
\right] $ is invertible, so the matrix $\left[
\begin{array}{cc}
j & 0 \\
0 & 1%
\end{array}
\right] $is a splittable monomorphism. So there exists a map $\left[
\begin{array}{cc}
x & y \\
z & w%
\end{array}
\right] $ that is a left inverse of $\left[
\begin{array}{cc}
j & 0 \\
0 & 1%
\end{array}
\right] $, and we obtain $xj=1$, so that $j$ is a splittable monomorphism
proving that $g_{i}$ is an irreducible morphism. We show in a similar way
that each $f_{i}$ is also an irreducible morphism. $\blacksquare \smallskip $

It turns out that each irreducible morphism in $gr_{0\Lambda }$ is either a
monomorphism or an epimorphism. For instance, assume that we have a map $%
f:E\rightarrow M$ that is irreducible in $gr_{0\Lambda }.$If $f$ is not
onto, then, since $X=Imf$ is again generated in degree zero, we have a
factorization $
\begin{array}{ccc}
\text{E} & \overset{\text{f}}{\longrightarrow } & \text{M} \\
\underset{\text{s}}{\searrow } &  & \underset{\text{t}}{\nearrow } \\
& \text{X} &
\end{array}
$and this implies that $s$ is a splittable monomorphism, hence $f$ is also a
monomorphism. This implies that if $0\rightarrow \sigma (M)\overset{f}{%
\rightarrow }E\overset{g}{\rightarrow }M\rightarrow 0$\allowbreak\ is an A-R
sequence in $K_{\Lambda },$ and $f=\left[ f_{1},f_{2},...f_{k}\right] $ and $%
g=$ $\left[ g_{1},g_{2},...g_{k}\right] ^{T}$, then for each $i$, the
morphisms $f_{i}$ and $g_{i}$ are either monomorphisms or
epimorphisms.\smallskip

We recall the following definition. Assume that the subcategory $\mathcal{C}$
has left Auslander-Reiten sequences, and let $M\in \mathcal{C}$. The cone of
$M$ is the subquiver of the A-R quiver of $\mathcal{C}$ , consisting of $M$
and its predecessors. Let now $M$ be a non projective Koszul module over an
indecomposable selfinjective Koszul algebra of Loewy length $r+1$, for
example the exterior algebra $\Lambda =\wedge ^{r+1}V,$where $V$ is a $K-$%
vector space of dimension $r+1$and $r\geq 2.$We will describe the cone of $M$
in $K_{\Lambda }$, by analyzing the cone of $M/J^{2}M$ in the A-R quiver of $%
gr_{\Lambda /J^{2}.}$ Namely, using the Auslander-Reiten quiver of $%
gr_{\Lambda /J^{2}}$ and $\left[ GMRSZ\right] $, we see that two of the
components of the A-R quiver of $gr_{0\Lambda /J^{2}}$ are of the form
\textsl{Z}$\Delta $, where $\Delta $ is the separated quiver of $\Lambda
/J^{2}$. For example in the case $\Lambda $ is the exterior algebra the
separated quiver of $\Lambda /J^{2}$ is the quiver $\bullet
\rightrightarrows \bullet $ with two vertices and $r+1$ arrows from the
first vertex to the second one. There are the ''preprojective'' component

\begin{center}
$
\begin{array}{cccccccc}
Y_{0}=\Lambda /J^{2} &  &  &  & Y_{2} &  &  &  \\
& \searrow \searrow &  & \nearrow \nearrow &  & \searrow \searrow &  &  \\
&  & Y_{1} &  &  &  & Y_{3} & ...%
\end{array}
$
\end{center}

and the ''preinjective'' component

\begin{center}
$
\begin{array}{cccccccc}
... & X_{3} &  &  &  & X_{1}=I &  &  \\
&  & \searrow \searrow &  & \nearrow \nearrow &  & \searrow \searrow &  \\
&  &  & X_{2} &  &  &  & X_{0}=K%
\end{array}
$
\end{center}

\begin{lemma}
Let $\Lambda $ be a Koszul algebra and let $f:M\rightarrow N$ be a map in $%
gr_{0\Lambda }$. Then the induced map $\overline{f}:M/J^{2}M\rightarrow
N/J^{2}N$ is a monomorphism if $f$ is a monomorphism and an epimorphism if $%
f $ is an epimorphism.
\end{lemma}

\textit{Proof. }The statement about epimorphisms is trivial. Assume $f$ is a
monomorphism, so we have an exact sequence $0\rightarrow M\rightarrow
N\rightarrow C\rightarrow 0$. It is easy to show that $J^{i}M=J^{i}N\cap M$
for all $i\geq 0$. The result follows from $\left[ GM1\right] $. $%
\blacksquare $

\begin{proposition}
Let $\Lambda $ be a selfinjective Koszul algebra, $M$ an indecomposable non
projective Koszul module, and let $0\rightarrow \sigma (M)\overset{f}{%
\rightarrow }E\overset{g}{\rightarrow }M\rightarrow 0$ be the A-R sequence
in $K_{\Lambda }$ ending at $M$, where the maps are $f=\left[
f_{1},f_{2},...f_{k}\right] $ and $g=$ $\left[ g_{1},g_{2},...g_{k}\right]
^{T}$. Then:

i) The induced sequence $0\rightarrow \sigma (M)\overset{\overline{f}}{%
\rightarrow }E/J^{2}E\overset{\overline{g}}{\rightarrow }M/J^{2}M\rightarrow
0$ is the A-R sequence ending at $M/J^{2}M$ in $gr_{0\Lambda /J^{2}}$.

ii) The number of indecomposable summands of E, equals the number of
indecomposable summands of $E/J^{2}E.$

iii) The irreducible morphisms $f_{i},$ $g_{i}$ are monomorphisms
(epimorphisms) if and only if the induced maps $\overline{f_{i}},$ $%
\overline{g}_{i}$ are monomorphisms (epimorphisms).
\end{proposition}

\textit{Proof. }i) Assume $X\rightarrow M/J^{2}M$ is not a nonsplittable
epimorphism in $gr_{0\Lambda /J^{2}}$. Using the equivalence mentioned
earlier, this map is induced by a homomorphism $h:Y\rightarrow M$ in $%
\mathcal{L}_{\Lambda }$ that is not a splittable epimorphism. Since $h$ can
be lifted to $E$, the result follows since the equivalence $\mathcal{L}%
_{\Lambda }\cong gr_{0\Lambda /J^{2}}$ restricted to $K_{\Lambda \text{,}}$
implies the indecomposability of $M/J^{2}M$.

ii) and iii) follow immediately from our previous remarks.$\blacksquare
\smallskip $

As an immediate application of this result we see that if $M$ is a Koszul
module of Loewy length 2, then the A-R sequence in $K_{\Lambda}$ ending at $%
M $coincides with the one in $gr_{0\Lambda/J^{2}}$. Let us assume now that $%
M $ is an indecomposable non projective module in $K_{\Lambda}.$ Then $%
M/J^{2}M $ is again an indecomposable $\Lambda/J^{2}-$ module, so by
analyzing the A-R quiver of $\Lambda/J^{2}$, we see that $M/J^{2}M$ either
belong to the preprojective or preinjective components of $%
gr_{0\Lambda/J^{2}}$ or is regular. We show that one of this possibilities
cannot exist:

\begin{lemma}
Let $\Lambda $ be an indecomposable selfinjective Koszul algebra of Loewy
length $r+1,$and $M$ an indecomposable non projective Koszul module. Then $%
M/J^{2}M$ can not lie in the preprojective component of $gr_{0\Lambda
/J^{2}} $.
\end{lemma}

\textit{Proof. }If it does, it follows from the previous result that $\sigma
(M)$ is also in the preprojective component.

This component has two kinds of $\tau -$orbits, one of them is the $\tau -$%
orbit of $P/J^{2}P$, with $P$ an indecomposable projective. If $\sigma (M)$
is in the orbit of $P/J^{2}P$, then by induction, $P/J^{2}P$ is the $\sigma $
of a Koszul module, therefore it must be a Koszul but we know that if $r>1$%
then $P/J^{2}P$ is not Koszul, and we obtain a contradiction. Assume $\sigma
(M)$ is in the orbit of $Y_{1}=\tau _{\Lambda /J^{2}}^{-1}(S)$ with S a
simple module generated in degree zero. This implies by induction $Y_{1}$ is
Koszul. On the other hand, it is easy to verify that over $\Lambda ,$ the
module $Y_{1}$ can not be a Koszul module since this would imply that $%
\sigma (Y_{1})$ is simple and this is impossible.$\blacksquare \smallskip $

We have enough to describe the shape of the Auslander Reiten quiver of $%
K_{\Lambda},$where $\Lambda$ is selfinjective Koszul. Using our previous
results and the preceding remarks, it is not hard to see that, if $M$ an
indecomposable non projective Koszul module, then the cones of $M$ in $%
K_{\Lambda}$ and of $M/J^{2}M$ in $gr_{0\Lambda/J^{2}}$ are isomorphic via
an isomorphism that takes monomorphisms into monomorphisms and epimorphisms
into epimorphisms. Moreover, since the $\Lambda-$modules $soc^{2}P=I,$where $%
P$ is an indecomposable projective injective and the simple $S$ are Koszul,
and the injective modules of $\Lambda/J^{2}$ are precisely the modules $I,$
the entire preinjective component of $gr_{0\Lambda/J^{2}}$ consists of
Koszul $\Lambda-$ modules. Thus the two kinds of $\tau-$orbits in this
component are in fact, $\sigma-$orbits. Putting together these facts gives
us the main result of this section:

\begin{theorem}
Let $\Lambda $ be an indecomposable selfinjective Koszul algebra of Loewy
length $r+1,$with $r>1.$The A-R quiver of $K_{\Lambda }$ has connected
components containing the indecomposable projective modules, a component
that coincides with the preinjective component of $gr_{0\Lambda /J^{2}},$and
all the remaining connected components are full subquivers of a quiver of
type \textsl{Z}$A_{\infty }.$
\end{theorem}

Since the preinjective component of $K_{\Lambda}$ consists only of Koszul
modules of Loewy length at most two, it turns out that each non projective
Koszul module of Loewy length three or higher lies in a ''regular''
component. The following result yields examples of modules in the regular
component of $K_{\Lambda}$ lying at the mouths of these components (see also
$\left[ GMRSZ\right] $, 3.2and 3.6).

\begin{proposition}
Let $\Lambda $ be as in the previous theorem, let $M$ be an indecomposable
non projective Koszul $\Lambda -$module and assume that $M$ has no
cogenerators in degree 1. Let $0\rightarrow \sigma (M)\overset{f}{%
\rightarrow }E\overset{g}{\rightarrow }M\rightarrow 0$ be the
Auslander-Reiten sequence in $K_{\Lambda }$ ending at $M$. Then, the middle
term $E$ is indecomposable.
\end{proposition}

\textit{Proof. }Assume that $E$ decomposes into $k>1$ indecomposable
summands. Then $f=\left[ f_{1},f_{2},...f_{k}\right] $ and $g=$ $\left[
g_{1},g_{2},...g_{k}\right] $ $^{T}$. It is clear that each composition $%
g_{i}f_{i}$ is nonzero. On the other hand, it was proved in $\left[ GMRSZ%
\right] $, that $\tau M$, and therefore $\sigma M$ too, are cogenereted in
degree $1$. Hence we obtain a contradiction to the fact that $M$ has no
cogenerators in degree $1$.$\blacksquare $

\section{Construction of locally free sheaves.}

In this section we will prove that for an indecomposable noetherian Artin
Shelter regular Koszul algebra $\Gamma$ of global dimension $r+1$ the
category of coherent sheaves $Qgr\Gamma$ has relative right almost split
sequences. Moreover, we will prove that the category of locally free sheaves
has relative right almost split sequences $\left[ AR\right] $. When we
specialize to $\Gamma=K\left[ x_{0},x_{1},...x_{r}\right] $ with $r\geq2$ we
will construct indecomposable vector bundles on \textsl{P}$_{r}$ of
arbitrary large rank.

\begin{lemma}
Let $\Gamma $ be a noetherian Artin Shelter regular Koszul algebra, $\pi
:gr\Gamma \rightarrow Qgr\Gamma $ the quotient functor and let $M$ be $%
\Gamma -$module with $\pi (M)$ indecomposable in $Qgr\Gamma .$ Then $%
End_{Qgr\Gamma }(\pi M)$ is local.
\end{lemma}

\textit{Proof. }Let $M$ be a graded torsion free module with $\pi (M)$
indecomposable. We know there exists an integer $k$ such that $M_{\geq k}%
\left[ k\right] =N$ is Koszul. Since $\pi (N\left[ -k\right] )=\pi (M),$ it
follows $N$ is indecomposable.

Hence; we may assume $M$ is Koszul up to shifting.

Let $L$ be a torsion free module we have:

\begin{center}
$Hom_{Qgr\Gamma }(\pi (M),\pi (L))_{0}=\underset{k\geq 0}{\cup }%
Hom_{gr\Gamma }(J^{k}M),L)_{0}.$
\end{center}

To see this, observe first that if $k<l,$ then $J^{l}M\subset J^{k}M$ and
the exact sequence:

$0\rightarrow J^{l}M\rightarrow J^{k}M$ $\rightarrow J^{k}M/$ $%
J^{l}M\rightarrow 0$ induces an exact sequence:

\begin{center}
$Hom_{gr\Gamma }(J^{k}M$ $/J^{l}M,L)_{0}\rightarrow Hom_{gr\Gamma
}(J^{k}M,L)_{0}$ $\rightarrow $ $Hom_{gr\Gamma }(J^{l}M,L)_{0}$,
\end{center}

with $Hom_{gr\Gamma }(J^{k}M$ $/J^{l}M,L)_{0}=0$.

It follows, $Hom_{gr\Gamma }(J^{k}M,L)_{0}$ $\subset Hom_{gr\Gamma
}(J^{l}M,L)_{0}$ and

\begin{center}
$\underrightarrow{\lim }Hom_{gr\Gamma }(J^{k}M,L)_{0}$ $=\underset{k\geq 0}{%
\cup }Hom_{gr\Gamma }(J^{k}M,L)_{0}.$
\end{center}

Assume $M$ is Koszul and torsion free.

Then $End_{Qgr\Gamma }(M)=\underset{k\geq 0}{\cup }Hom_{gr\Gamma
}(J^{k}M,M)_{0}.$

Let $f\in Hom_{gr\Gamma }(J^{k}M,M)_{0}.$ Then $\mathit{Im}f\subset J^{k}M.$
Set $f$ $^{\prime }$ to be the restriction of $f$ to $J^{k}M.$

Let $\Lambda $ be the Yoneda algebra of $\Gamma $ and let $F:K_{\Lambda
^{op}}\rightarrow K_{\Gamma }$ be the Koszul duality.

We have natural isomorphisms:

\begin{center}
$Hom_{\Lambda ^{op}}(\Omega ^{k}F^{-1}(M),\Omega ^{k}F^{-1}(M))_{0}\cong
Hom_{\Gamma }(J^{k}M,J^{k}M)_{0}\cong $

$Hom_{\Lambda ^{op}}(F^{-1}(M),F^{-1}(M))_{0}\cong Hom_{\Gamma }(M,M)_{0}.$
\end{center}

It follows, $End_{\Gamma }(J^{k}M)_{0}$ is local and $f^{\prime }$ is either
nilpotent or invertible.

Then $f$ is either nilpotent or invertible.$\blacksquare $

\begin{proposition}
Let $\Gamma $ be an indecomposable noetherian Artin Shelter regular Koszul
algebra of global dimension $r+1,$ with Yoneda algebra $\Lambda .$ Let $%
F:K_{\Lambda }\rightarrow K_{\Gamma ^{op}}$ be the Koszul duality. Then the
category of non simple Koszul $\Gamma ^{op}-\mathit{mod}$ules has right almost
split sequences.
\end{proposition}

\textit{Proof. }Let $M$ be an indecomposable non simple Koszul $\Gamma
^{op}- $ module. Then there exists a non projective Koszul $\Lambda -$
module $N$ such that $F(N)=M.$

Let $0\rightarrow \tau (N)\rightarrow E\rightarrow N\rightarrow 0$ be the
almost split sequence in $gr\Lambda .$ The module $\tau (N)$ is generated in
degree $-r+1$ and $\tau (N)\left[ -r+1\right] $ is Koszul.

It was proved in $\left[ GMRSZ\right] ,$ $J^{r-1}\tau (N)$ is Koszul $%
0\rightarrow \tau (N)_{\geq 0}\rightarrow E_{\geq 0}\rightarrow N\rightarrow
0$ is an almost split sequence in $K_{\Lambda }$ and $\tau (N)_{\geq
0}=J^{r-1}\tau (N)$ is indecomposable.

Applying Koszul duality $F$ we obtain an almost split sequence in $K_{\Gamma
^{op}}$:

$0\rightarrow F(N)\rightarrow F($ $E_{\geq 0})\rightarrow $ $F(J^{r-1}\tau
(N))\rightarrow 0$ .$\blacksquare \smallskip \smallskip $

Let $\Gamma$ be a noetherian Artin Schelter regular Koszul algebra and $%
\pi:gr_{\Gamma^{op}}\rightarrow Qgr_{\Gamma^{op}}$ the quotient functor.
Denote by $\overset{\wedge}{K}_{\Gamma^{op}}\left[ n\right] $ the
subcategory of $Qgr_{\Gamma^{op}}$ defined as: $\overset{\wedge}{K}%
_{\Gamma^{op}}\left[ n\right] =\left\{ \pi(M\left[ n\right] )\mid M\text{ is
a Koszul }\Gamma^{op}\text{-}\mathit{mod}\text{ule}\right\} $ it follows by
Proposition 5, the category $Qgr_{\Gamma^{op}}$ is equal to $\underset{n\in Z%
}{\cup}$ $\overset{\wedge}{K}_{\Gamma^{op}}\left[ n\right] .$ With this
notation we have the following:

\begin{definition}
An object $\overset{\thicksim }{M}$ in $Qgr_{\Gamma ^{op}}$ has a relative
right almost split sequence if given an integer $n$ such that $\overset{%
\thicksim }{M}\in $ $\overset{\wedge }{K}_{\Gamma ^{op}}\left[ n\right] $ ,
then there exist $\overset{\thicksim }{L},\overset{\thicksim }{N\in }\overset%
{\wedge }{K}_{\Gamma ^{op}}\left[ n\right] $ and a short exact sequence: $%
0\rightarrow \overset{\thicksim }{M}\rightarrow \overset{\thicksim }{L}%
\rightarrow \overset{\thicksim }{N}\rightarrow 0$ which is an almost split
sequence in $\overset{\wedge }{K}_{\Gamma ^{op}}\left[ n\right] .$
\end{definition}

\begin{theorem}
Let $\Gamma $ be an indecomposable noetherian Artin Shelter regular Koszul
algebra. Then the category of sheaves $Qgr\Gamma ^{op}$ has relative right
almost split sequences.
\end{theorem}

\textit{Proof. }Let $\overset{\thicksim }{M}$ be an indecomposable object in
$Qgr\Gamma ^{op}$ and $\pi :gr\Gamma ^{op}\rightarrow Qgr\Gamma ^{op}$ be
the quotient functor. Then there exists an indecomposable torsion free $%
\Gamma -$module $M$ such that $\pi (M)=\overset{\thicksim }{M}.$

Since for some integer $k$ the truncated module $M_{\geq k}\left[ k\right]
=X $ is Koszul and $\pi (X\left[ -k\right] )=\pi (M),$ hence; $X$ is
indecomposable..

By proposition 10, there exists an almost split sequence in $K_{\Gamma
^{op}} $ of the form: $0\rightarrow X\overset{i}{\rightarrow }Y\overset{p}{%
\rightarrow }Z\rightarrow 0$ and an exact sequence:

\begin{center}
$0\rightarrow X\left[ -k\right] \rightarrow Y\left[ -k\right] \rightarrow Z%
\left[ -k\right] \rightarrow 0$
\end{center}

The functor $\pi $ is exact $\left[ P\right] ,\left[ S2\right] .$

Hence; there exists an exact sequence:

\begin{center}
$0\rightarrow \pi (X\left[ -k\right] )\overset{\pi (i)}{\rightarrow }\pi (Y%
\left[ -k\right] )\overset{\pi (p)}{\rightarrow }\pi (Z\left[ -k\right]
)\rightarrow 0.\smallskip $
\end{center}

Assume there exits $h:\pi (Y\left[ -k\right] )\rightarrow \pi (X\left[ -k%
\right] )$ such that $h\pi (i)=1$, where $h\in \cup Hom(J^{t}Y\left[ -k%
\right] ,X\left[ -k\right] )).$ Then the\smallskip\ sequence:\smallskip
\smallskip\ $0\rightarrow J^{t}X\left[ -k\right] \rightarrow J^{t}Y\left[ -k%
\right] \rightarrow J^{t}Z\left[ -k\right] \rightarrow 0$ splits: It
follows\smallskip \newline
the sequence: $0\rightarrow \Omega ^{t}F^{-1}(Z)\left[ -k\right] \rightarrow
\Omega ^{t}F^{-1}(Y)\left[ -k\right] \rightarrow \Omega ^{t}F^{-1}(X)\left[
-k\right] \rightarrow 0$ splits. Then the sequence:\smallskip\ $0\rightarrow
F^{-1}(Z)\left[ -k\right] \rightarrow F^{-1}(Y)\left[ -k\right] \rightarrow
F^{-1}(X)\left[ -k\right] \rightarrow 0$ splits. It follows\smallskip\ $%
0\rightarrow X\left[ -k\right] \rightarrow Y\left[ -k\right] \rightarrow Z%
\left[ -k\right] \rightarrow 0$ splits.

We must prove the sequence is almost split in $\overset{\wedge }{K}_{\Gamma
^{op}}\left[ -k\right] $.

Let $h:\pi (W\left[ -k\right] )\rightarrow \pi (Z\left[ -k\right] )$ be a
non splittable map in $\overset{\wedge }{K}_{\Gamma ^{op}}\left[ -k\right] .$
The map\smallskip \newline
$h$ belongs to $\cup $ $Hom(J^{t}W\left[ -k\right] ,Z\left[ -k\right] ))$,
then $h:J^{t}W\left[ -k\right] \rightarrow J^{t}Z\left[ -k\right] $ does not
split.

Applying the functor $F^{-1}$we obtain the following exact diagram:\medskip

$
\begin{array}{ccccccc}
0\rightarrow & \Omega ^{t}F^{-1}(Z)\left[ -k\right] & \rightarrow & \Omega
^{t}F^{-1}(Y)\left[ -k\right] & \rightarrow & \Omega ^{t}F^{-1}(X) \left[ -k%
\right] & \rightarrow 0 \\
& \downarrow &  &  &  &  &  \\
& \Omega ^{t}F^{-1}(W)\left[ -k\right] &  &  &  &  &
\end{array}
\medskip $

Applying $\Omega ^{-t}$ to the diagram and using the fact $0\rightarrow
F^{-1}(Z)\left[ -k\right] \rightarrow F^{-1}(Y)\left[ -k\right] $

$\rightarrow F^{-1}(X)\left[ -k\right] \rightarrow 0$ is almost split in $%
K_{\Lambda }\left[ -k\right] $ we obtain an exact diagram:\medskip

$
\begin{array}{ccccccc}
0\rightarrow & F^{-1}(Z)\left[ -k\right] & \rightarrow & F^{-1}(Y)\left[ -k%
\right] & \rightarrow & F^{-1}(X)\left[ -k\right] & \rightarrow 0 \\
& \downarrow & \swarrow &  &  &  &  \\
& F^{-1}(W)\left[ -k\right] &  &  &  &  &
\end{array}
\medskip $

Applying $\Omega ^{t}$ to the diagram we obtain and extension of the map
\smallskip

\begin{center}
$\Omega ^{t}F^{-1}$(Z)$\left[ -k\right] \rightarrow \Omega ^{t}F^{-1}(W)%
\left[ -k\right] .\smallskip $
\end{center}

Applying the functor $F$ we obtain a lifting of the map $h:J^{t}W\left[ -k%
\right] \rightarrow J^{t}Z\left[ -k\right] $ to $g:J^{t}W\left[ -k\right]
\rightarrow Y\left[ -k\right] $. Finally, applying the functor $\pi $ to the
diagram we obtain\smallskip \newline
a lifting $\pi (g):\pi (W\left[ -k\right] )\rightarrow \pi (Y\left[ -k\right]
)$ of $h$ as claimed.

We have proved $0\rightarrow \pi (M)\overset{\pi (i)}{\rightarrow }\pi (L)%
\overset{\pi (p)}{\rightarrow }\pi (N)\rightarrow 0$ is a relative almost
split sequence in $Qgr\Gamma ^{op}.\blacksquare $

\begin{corollary}
Let $Coh(sl{P}_{r})$ be the category of coherent sheaves on projective
space $sl{P}_{r}.$ Then $Coh(sl{P}_{r})$ has relative right
almost split sequences.$\blacksquare $
\end{corollary}

\bigskip Let $\Gamma$ be an Artin Schelter regular Koszul algebra, $L$ and $%
M $ finitely generated Koszul $\Gamma^{op}-$module such that $\pi(L\left[ n%
\right] )\cong\pi(M\left[ m\right] )$, and assume $n\geq m.$ Then there
exits an integer $k$ such that $L\left[ n\right] _{\geq k}=L_{\geq n+k}\left[
n\right] \cong M\left[ m\right] _{\geq k}=M_{m+k}\left[ m\right] $ or $%
J^{n+k}L\left[ n+k\right] \cong J^{m+k}M\left[ m+k\right] .$

Applying the duality $F^{-1}$ we obtain: $\Omega^{n+k}F^{-1}L\left[ n+k%
\right] \cong\Omega^{m+k}F^{-1}M\left[ m+k\right] .$

Therefore: $\Omega^{n-m}F^{-1}L\left[ n-m\right] \cong F^{-1}M.$ Applying $F$
it follows: $J^{n-m}L\left[ n-m\right] $

$\cong M.$

The following leemma will be needed in the next proposition:

\begin{lemma}
Let $\Lambda $ be an indecomposable graded selfinjective quiver algebra of
Loewy length $r+1$. Then for any non negative integer $k$ we have an
isomorphism: $(D(\Omega ^{k}\Lambda _{0}))^{\ast }\left[ r+1\right] \cong
\Omega ^{k}\Lambda _{0\text{.}}$
\end{lemma}

\textit{Proof. }Let $S$ be a graded simple generated in degree zero. Then By
$\left[ M3\right] ,$ $D(S)$ is also generated in degree zero. Let $P$ be the
projective cover of $D(S)$ applying the functor $(-)^{\star }=Hom_{\Lambda
}(-,\Lambda )$ to the exact sequence: $P\rightarrow D(S)\rightarrow 0$ we
obtain an exact sequence: $0\rightarrow (D(S))^{\star }\rightarrow P^{\star
} $ with $P^{\star }$ an indecomposable projective generated in degree zero
of Loewy length $r+1$. Therefore $(D(S))^{\star }\left[ r+1\right] $ is a
simple module generated in degree zero. The Nakayama equivalence: $%
S\rightarrow (D(S))^{\star }$ $\left[ r+1\right] $ induces a bijection of
the graded simple modules generated in degree zero, hence $\Lambda _{0}\cong
(D(\Lambda _{0}))^{\star }\left[ r+1\right] .$

Applying $(D(-))^{\star }\left[ r+1\right] $ to the minimal projective
resolution of $\Lambda _{0}$ we obtain a minimal projective resolution of $%
(D(\Lambda _{0}))^{\star }\left[ r+1\right] $. It follows: $(D(\Omega
^{k}\Lambda _{0}))^{\ast }\left[ r+1\right] \cong \Omega ^{k}\Lambda _{0}$.$%
\blacksquare \smallskip $

\begin{theorem}
Let $\Gamma $ be an indecomposable noetherian Artin Shelter regular Koszul
algebra of global dimension $r+1$. Then the category of locally free sheaves
of finite rank has relative right almost split sequences.
\end{theorem}

\textit{Proof. }Let $\Lambda $ be the Yoneda algebra of $\Gamma $ and $%
F:K_{\Lambda }\rightarrow K_{\Gamma ^{op}}$ Koszul duality. Let $\overset{%
\thicksim }{M}\in Qgr\Gamma ^{op}$ be an indecomposable locally free sheaf,
we may assume $\overset{\thicksim }{M}=\pi (M)$, where $\pi :gr\Gamma
^{op}\rightarrow Qgr\Gamma ^{op}$ is the quotient functor and $M$ is an
indecomposable module.

Then there exists an indecomposable non projective Koszul module $X$ and
some integer $n$ with $F(X\left[ n\right] )=M.$

Let $0\rightarrow \tau (X)\rightarrow E\rightarrow X\rightarrow 0$ be the
almost split sequence in $gr\Lambda $ and $0\rightarrow J^{r-1}\tau
(X)\rightarrow E_{\geq 0}\rightarrow X\rightarrow 0$ the almost split
sequence in $K_{\Lambda }.$

Since $\pi (F(X))$ is locally free $Ext_{\Lambda }^{j}(\Omega ^{s}\Lambda
_{0}\left[ s\right] ,X)_{0}=0$ for $0<j\leq r$ and $s>>0.$

There is a chain of isomorphisms:\smallskip

$Hom_{\Lambda }(\Omega ^{s+j}\Lambda _{0}\left[ s\right] ,$ $D(X^{\ast })%
\left[ -r-1\right] )_{0}\cong Hom_{\Lambda }(\Omega ^{s+j}\Lambda _{0},$ $%
D(X^{\ast })\left[ -r-s-1\right] )_{0}\cong \smallskip $

$Hom_{\Lambda }((D(\Omega ^{s+j}\Lambda _{0}))^{\ast },$ X$\left[ -r-s-1%
\right] )_{0}\cong Hom_{\Lambda }((D(\Omega ^{s+j}\Lambda _{0}))^{\ast }%
\left[ r+1\right] ,$ $X\left[ -s\right] )_{0}\smallskip $

$\cong Hom_{\Lambda }(\Omega ^{s+j}\Lambda _{0},$ $X\left[ -s\right]
)_{0}\cong Hom_{\Lambda }(\Omega ^{s+j}\Lambda _{0}\left[ s\right] ,$ $%
X)_{0}=0$ for $s>>0.\smallskip $

But $\tau (X)\cong \Omega ^{2}D(X^{\ast }),$ hence; $\tau (X)\left[ -r+1%
\right] \cong \Omega ^{2}D(X^{\ast })\left[ -r-1\right] \left[ 2\right]
\smallskip $

It follows by Corollary 2 $\tau (X)\left[ -r+1\right] $ is such that $\pi
(F(\tau (X)\left[ -r+1\right] ))$ is locally free.

By Proposition 5 $J^{r-1}\tau (X)$ is such that $\pi (F(J^{r-1}\tau (X)))$
is locally free and by Proposition 4 $\pi (F(E_{\geq 0}))$ is locally free.

It follows $0\rightarrow \pi (F(X))\rightarrow \pi (F(E_{\geq
0}))\rightarrow \pi (F(J^{r-1}\tau (X)))\rightarrow 0$ is a relative almost
split sequence of locally free sheaves, hence;. $0\rightarrow \pi (F(X)\left[
n\right] )\rightarrow \pi (F(E_{\geq 0})\left[ n\right] )\rightarrow \pi
(F(J^{r-1}\tau (X))\left[ n\right] )\rightarrow 0$ is relative almost split.$%
\blacksquare $

\begin{corollary}
The category of vector bundles \textit{V(}\textsl{P}$_{r})$ on projective
space \textsl{P}$_{r}$ has relative right almost split sequences.$%
\blacksquare $
\end{corollary}

Let $\Gamma $ and $\Lambda $ be like in the theorem. The almost split
sequences in $Qgr\Gamma ^{op}$ are related with the almost split sequences
in $Gr_{\Lambda /J^{2}}$, using the results of Section 3 we have the
following version of the main theorem of Section 3:.

\begin{theorem}
Let $\Gamma $ be an indecomposable noetherian Artin Schelter regular Koszul
algebra of global dimension $r$ with $r>2,$with Yoneda algebra $\Lambda $,
let $F:K_{\Lambda }\rightarrow K_{\Gamma ^{op}}$ be the Koszul duality and
let $M$ be an indecomposable Koszul $\Lambda -$ module. Then the relative
Auslander Reiten component of $\pi F(M)$ is either contained in the
nonegative part of a quiver of type \textsl{Z}$\Delta ,$where $\Delta $ is
the separated quiver of $\Gamma /J^{2}$ or of the form:
\end{theorem}

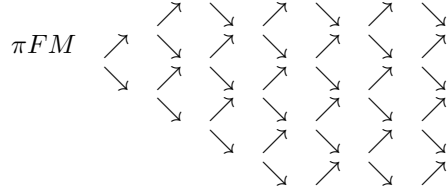
\begin{center}
$
\begin{array}{cccccccc}
&  & \nearrow & \searrow & \nearrow & \searrow & \nearrow & \searrow \\
\pi FM & \nearrow & \searrow & \nearrow & \searrow & \nearrow & \searrow &
\nearrow \\
& \searrow & \nearrow & \searrow & \nearrow & \searrow & \nearrow & \searrow
\\
&  & \searrow & \nearrow & \searrow & \nearrow & \searrow & \nearrow \\
&  &  & \searrow & \nearrow & \searrow & \nearrow & \searrow \\
&  &  &  & \searrow & \nearrow & \searrow & \nearrow%
\end{array}
$
\end{center}

\textit{If }$\pi (F(M))$\textit{\ is locally free, then all the sheaves in
the diagram are locally free.}

\textit{For all the modules }$M$\textit{\ in the nonegative part of a quiver
of type Z}$\Delta ,$\textit{where }$\Delta $\textit{\ is the separated
quiver of }$\Gamma /J^{2}$\textit{\ the sheaf }$\pi (F(M))$\textit{\ is
locally free.}$\blacksquare $

As a corollary we have:

\begin{theorem}
Let $\Gamma =K\left[ x_{0},x_{1},...x_{r}\right] $ be the polynomial algebra
with $r>1$ and $\Lambda =K<x_{0},x_{1},...x_{r}>/(x_{i}^{2},$ $%
x_{i}x_{j}+x_{j}x_{i})$ the exterior algebra, let $F:K_{\Lambda }\rightarrow
K_{\Gamma ^{op}}$ be the Koszul duality and let $M$ be an indecomposable
Koszul $\Lambda -$ module. Then the relative Auslander Reiten component of $%
\pi F(M)$ is either contained in the non negative part of the quiver \textsl{%
Z}$\Delta ,$ where $\Delta $ is the quiver: $\bullet \rightrightarrows
\bullet $ with $r+1$ arrows, or of the form:
\end{theorem}

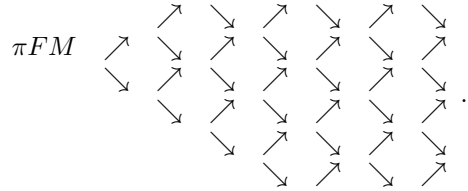
\begin{center}
$
\begin{array}{cccccccc}
&  & \nearrow & \searrow & \nearrow & \searrow & \nearrow & \searrow \\
\pi FM & \nearrow & \searrow & \nearrow & \searrow & \nearrow & \searrow &
\nearrow \\
& \searrow & \nearrow & \searrow & \nearrow & \searrow & \nearrow & \searrow
\\
&  & \searrow & \nearrow & \searrow & \nearrow & \searrow & \nearrow \\
&  &  & \searrow & \nearrow & \searrow & \nearrow & \searrow \\
&  &  &  & \searrow & \nearrow & \searrow & \nearrow%
\end{array}
.$
\end{center}

\textit{If }$\pi (F(M))$\textit{\ is locally free, then all the sheaves in
the diagram are locally free.}

\textit{If }$M$\textit{\ is contained in the non negative part of the quiver
Z}$\Delta ,$\textit{\ where }$\Delta $\textit{\ is the quiver: }$\bullet
\rightrightarrows \bullet $\textit{\ with }$r+1$\textit{\ arrows, then }$\pi
(F(M))$\textit{\ is locally free.}$\blacksquare $

\begin{corollary}
There are indecomposable vector bundles on \textsl{P}$_{r}$ , with $r>1,$ of
arbitrary high rank.
\end{corollary}

\textit{Proof. }If we have an exact sequence of vector bundles: $%
0\rightarrow X\rightarrow Y\rightarrow Z\rightarrow 0$ it is clear that $rk$
$(Y)=rk$ $(X)+rk$ $(Z),$ where $rk$ $(Y)$ dentes the rank of $Y.$

Let $M$ be an indecomposable Koszul $\Lambda -$ module and assume $%
Hom_{\Lambda }(\Lambda _{0}\left[ -1\right] ,M)_{0}=0,$ for example $M=J%
\left[ 1\right] .$ Then the relative Auslander Reiten component of $\pi F(M)$
is of the form:

\begin{center}
$
\begin{array}{ccccccccc}
\pi F(M) &  &  &  & \pi F(\sigma M) &  &  &  & \pi F(\sigma ^{2}M) \\
& \searrow &  & \nearrow &  & \searrow &  & \nearrow &  \\
&  & \pi F(M_{1}) &  &  &  & \pi F(\sigma M_{1}) &  &  \\
&  &  & \searrow &  & \nearrow &  & \searrow &  \\
&  &  &  & \pi F(M_{2}) &  &  &  & \pi F(\sigma M_{2}) \\
&  &  &  &  & \searrow &  & \nearrow &  \\
&  &  &  &  &  & \pi F(M_{3}) &  &
\end{array}
$
\end{center}

Set $M=M_{0}.$

For any $j\geq 0$ the following equality hold:

\begin{center}
$rk(\pi F(\sigma ^{j}M_{1}))=rk(\pi F(\sigma ^{j+1}M))+rk(\pi F(\sigma
^{j}M))$.
\end{center}

Assume for all $j\geq 0$ there is an equality:

\begin{center}
$rk(\pi F(\sigma ^{j}M_{i}))=rk(\pi F(\sigma ^{j}M_{i-1}))$ $+rk(\pi
F(\sigma ^{i+j}M))$
\end{center}

We have the following commutative diagram:

\begin{center}
$
\begin{array}{ccccc}
&  & \pi F(\sigma ^{j+1}M_{i-1}) &  &  \\
& \nearrow &  & \searrow &  \\
\pi F(\sigma ^{j}M_{i}) &  &  &  & \pi F(\sigma ^{j+1}M_{i}) \\
& \searrow &  & \nearrow &  \\
&  & \pi F(\sigma ^{j}M_{i+1}) &  &
\end{array}
$
\end{center}

Where the sum of the two terms in the middle is the middle term of an exact
sequence. It follows:

\begin{center}
$rk(\pi F(\sigma ^{j}M_{i}))+rk(\pi F(\sigma ^{j+1}M_{i}))=rk(\pi F(\sigma
^{j}M_{i+1}))+rk(\pi F(\sigma ^{j+1}M_{i-1}))$
\end{center}

Hence; $rk(\pi F(\sigma ^{j}M_{i+1}))$=$rk(\pi F(\sigma ^{j}M_{i}))+rk(\pi
F(\sigma ^{j+1}M_{i}))-rk(\pi F(\sigma ^{j+1}M_{i-1}))$

\begin{center}
$rk(\pi F(\sigma ^{j}M_{i+1}))=rk(\pi F(\sigma ^{j}M_{i}))+rk(\pi F(\sigma
^{i+j+1}M))$
\end{center}

It follows: $rk(\pi F(\sigma ^{j}M_{i}))<rk(\pi F(\sigma ^{j}M_{i+1}))$ for
all $j\geq 0$ and $i\geq 0.\blacksquare $

\textit{\
\bibliographystyle{AABBRV}
\bibliography{acompat}
}

\begin{center}
REFERENCES
\end{center}

$\left[ ADL\right] $ Agoston, I., Dlab, V., Luk\'{a}s, E. Homological
duality and quasi hereditary, Can. J. Math. 48 (1996), 897-917.

$\left[ AS\right] $ Artin, M.; Shelter W.; Graded algebras of dimension 3,
Advances in Mathematics, 66 (1987), 172-216.

$\left[ ABPRS\right] $ Auslander, M.; Bautista. R.; Platzeck, M.I.; Reiten,
I.; Smalo, S., Canadian J. Math. 31 (1979), No. 5, 942-960.

$\left[ AR\right] $ Auslander, M.; Reiten, I., Representation Theory of
Artin algebras III, Almost split sequences, Comm. in Algebra 3, (1975),
239-294.

$\left[ BGK\right] $ Baranovsky, V.; Ginzburg, V.; Kuznetsov, A.; Quiver
varieties and a noncommutative \textsl{P}$_{2},$ preprint (2001).

$\left[ Be\right] $ Belinson, A. Coherent Sheaves on \textsl{P}$_{n}$ and
problems of linear algebra. Funkts. Anal. Prilozh. 12, No. 3 (1978),
English. Trans: Funct. Anal. Appl. 12 (1979), 214-216.

$\left[ BGS\right] $ Belinson, A.; Ginzburg, V.; Soergel, W., Koszul duality
patterns in representation theory, J. Amer. Math. Society, 9, no. 2.,(1996),
473-527.

$\left[ Bo\right] $ Bondal, A., I., Representation of Associative Algebras
and Coherent Sheaves. Izv. Akad. Nauk SSSR, Ser. Mat. 53, No. 1 (1989),
25-44. English translation: Math USSR, Izv. 34 (1990), 23-42.

$\left[ DM\right] $ Dowbor, P.; Meltzer, H., On equivalences of
Bernstain-Gelfand-Gelfand, Beilinson and Happel. Comm. in Algebra 20 (9),
(1992), 2513-2532.

$\left[ G\right] $ Gelfand, S., I., Sheaves on \textsl{P}$_{n}$ and problems
of linear algebra. Appendix, In. C. Okonek, M. Schneider and H. Spindler,
Vector bundles on complex projective spaces, Moscow, Mir (1984), 278-305.

$\left[ GMRSZ\right] $ Green, E.,L.;\ Mart\'{\i}nez-Villa, R.; Reiten, I.;
Solberg, \O .; Zacharia, D., On modules with linear presentations, J.
Algebra 205, (1998), no- 2, 578-604.

$\left[ GM1\right] $ Green, E., L.; Mart\'{\i}nez-Villa, R., Koszul and
Yoneda algebras, Representation Theory of Algebras, 247-297, CMS Conf.
Proc., 18, Amer. Math. Soc. (1996).

$\left[ GM2\right] $ Green, E., L.; Mart\'{\i}nez-Villa, R., Koszul and
Yoneda algebras II, Algebras and Modules, II, 227-244, CMS Conf. Proc., 24,
Amer. Math. Soc. (1998).

$\left[ GMT\right] $ Guo, J., Y.; Mart\'{\i}nez-Villa, R.; Takane, M.,
Koszul generalized Auslander regular algebras, Algebras and Modules, II,
263-283, CMS Conf. Proc., 24, Amer. Math. Soc. (1998).

$\left[ H1\right] $ Hartshorne, R., Algebraic Vector Bundles on Projective
Space: A Problem List, Topology Vol. 18, 117-128, (1979).

$\left[ H2\right] $ Hartshorne, R., Algebraic Geometry, Springer, Graduate
Texts in Math. 52, (1997).

$\left[ L\right] $ Le Potier, J. Lectures on Vector Bundles, Cambridge
studies in advanced mathematics, 54, (1997)

$\left[ M1\right] $ Mart\'{\i}nez-Villa, R., Applications of Koszul
algebras: the preprojective algebra, Representation Theory of Algebras,
487-504, CMS Conf. Proc., 18, Amer. Math. Soc. (1996).

$\left[ M2\right] $ Mart\'{\i}nez-Villa, R., Serre Duality for Generalized
Auslander Regular Algebras, Contemporary Math. 229, (1998) 237-263.

$\left[ M3\right] $ Mart\'{\i}nez-Villa, R., Graded, Selfinjective, and
Koszul algebras, J. Algebra 215 (1999), no.1, 34-72.

$\left[ M4\right] $ Mart\'{\i}nez-Villa, R., Skew group algebras and their
Yoneda algebras, J. Okayama U. Vol. 43, 1-10, (2001)

$\left[ M5\right] $ Mart\'{\i}nez-Villa, R. Koszul algebras and the
Gorenstein condition, Representation of Algebras, Lecture Notes in Pure and
Applied Math. Math. Vol. 224, Dekker, (2001)

$\left[ MM\right] $ Mart\'{\i}nez-Villa, R.; Martsinkovsky, A., Sheaf
Cohomology and Tate Cohomology for Koszul algebras, preprint, (2001).

$\left[ MZ\right] $ Mart\'{\i}nez-Villa, R.; Zacharia, D.; Approximations
with modules having linear resolutions, preprint, (2001).

$\left[ P\right] $ Popescu, N., Abelian Categories with Applications to Ring
and Modules, Academic Press (1973).

$\left[ \Pr\right] $ Priddy, S., Koszul resolutions, Trans. Amer. Math. Soc.
152 (1970), 39-60.

$\left[ RR\right] $ Reiten, I..; Riedtmann, C., Skew Group Algebras in the
Representation Theory of Artin Algebras, J. of Algebra 92, 224-282, (1985)

$\left[ R\right] $ Rudakov, A., Helices and Vector Bundles, LMS Lect. Notes
Series 148, Cambridge, 1990.

$\left[ Sc\right] $ Schneider, M., Holomorphic Vector Bundles on \textsl{P}$%
_{n}$ , S\'{e}minaire Bourbaki, 31 ann\'{e}s, 1978/1979, no. 530.

$\left[ Se\right] $ Serre, J., P., Faisceaux alg\'{e}briques coh\'{e}rents,
Ann. Math. 61 (1955), 197-278.

$\left[ S1\right] $ Smith, P., Some finite-dimensional algebras related to
elliptic curves, Representation Theory of Algebras and related topics,
315-348, CMS Conf. Proc., 19, Amer. Math. Soc. (1996).

$\left[ S2\right] $ Smith, P., Non- Commutative Algebraic Geometry,
preprint, University of Washington.

\end{document}